\documentclass[12 pt]{amsart}

\usepackage{amscd,amsmath,latexsym,amssymb,amscd}
\usepackage[mathscr]{euscript}
\usepackage[all]{xy}
\usepackage{layout}
\usepackage{pb-diagram}
\usepackage{pstricks}
\usepackage{enumerate}
\usepackage{dsfont}

\newtheorem{theorem}{Theorem}[section]
\newtheorem{corollary}[theorem]{Corollary}
\newtheorem{lemma}[theorem]{Lemma}
\newtheorem{proposition}[theorem]{Proposition}
\theoremstyle{definition}
\newtheorem{definition}[theorem]{Definition}

\theoremstyle{definition}
\newtheorem{example}[theorem]{Example}

\def\Proof{\medskip\noindent{\bf Proof: }}

\def\Z{\mathbb{Z}}

\def\C{\mathbb{C}}
\def\D{\mathbb{D}}

\def\Q{\mathbb{Q}}
\def\R{\mathbb{R}}

\def\S{\mathbb{S}}
\newcommand{\BONE}{\mathds 1}

\def\Ca{\mathcal{C}}

\def\I{\mathcal{I}}

\def\Fi{\mathcal{F}}
\def\H{\mathcal{H}}

\def\W{\mathcal{W}}
\def\O{\mathcal{O}}
\def\P{\mathcal{P}}
\def\K{\mathcal{K}}
\def\lR{\mathcal{R}}
\def\lZ{\mathcal{Z}}
\def\bK{\mathbb{K}}

\def\Mod{\mathcal{M}od}

\def\g{\mathfrak{g}}

\def\t{\mathfrak{t}}
\def\gl{\mathfrak{gl}}
\def\su{\mathfrak{su}}
\def\ux{\underline{x}}

\DeclareMathOperator{\Hom}{Hom}

\begin{document}

\title[Equivariant K-theory of group actions with maximal rank isotropy]
{Equivariant K-theory of compact Lie group actions with
maximal rank isotropy}

\author[A.~Adem]{Alejandro Adem}
\address{Department of Mathematics,
University of British Columbia, Vancouver BC V6T 1Z2, Canada}
\email{adem@math.ubc.ca}

\author[J.~M.~G\'omez]{Jos\'e Manuel G\'omez}
\address{Department of Mathematics,
Johns Hopkins University, Baltimore, MD 21218, USA}
\email{jgomez@math.jhu.edu}

\date{\today}
\begin{abstract}
Let $G$ denote a compact connected Lie group with torsion--free
fundamental group acting on a compact space $X$ such that all the
isotropy subgroups are connected subgroups of maximal rank. Let
$T\subset G$ be a maximal torus with Weyl group $W$. If the
fixed--point set $X^T$ has the homotopy type of a finite $W$--CW
complex, we prove that the rationalized complex equivariant
$K$--theory of $X$ is a free module over the representation ring of
$G$. Given additional conditions on the $W$--action on the
fixed-point set $X^T$ we show that the equivariant $K$--theory of
$X$ is free over $R(G)$. We use this to provide computations for a
number of examples, including the ordered $n$--tuples of commuting
elements in $G$ with the conjugation action.
\end{abstract}

\maketitle

\section{Introduction}

Let $G$ denote a compact connected Lie group with torsion--free
fundamental group. Suppose that $G$ acts on a compact space $X$ so
that each isotropy subgroup is a connected subgroup of maximal rank.
In this article we study the problem of computing $K^*_G(X)$, the
complex $G$--equivariant $K$--theory of $X$. Our work is motivated
by the examples given by spaces of ordered commuting $n$-tuples in
compact matrix groups such as $SU(r)$, $U(q)$ and $Sp(k)$ with the
conjugation action.

\medskip

Given such a $G$--space $X$, if $T\subset G$ is a maximal torus,
then $N_{G}(T)$ acts on the $T$--fixed points $X^{T}$. This in turn
yields a corresponding action of the Weyl group $W= N_G(T)/T$ on the
fixed--point set. It can be shown (Theorem \ref{CW structures}) that
if $X^{T}$ has the homotopy type of a finite $W$--CW complex, then
this determines a $G$--CW complex structure on $X$, and its
$G$--equivariant $K$--theory becomes a tractable invariant. The
following theorem provides a calculation for $K^*_G(X)\otimes\mathbb
Q$.

\begin{theorem}\label{mainintroduction1}
Let $G$ be a compact connected Lie group with torsion--free
fundamental group and $T\subset G$ a maximal torus. Suppose that $X$
is a compact $G$--CW complex with connected maximal rank isotropy
subgroups, then $K^{*}_{G}(X)\otimes \Q$ is a free module over
$R(G)\otimes \Q$ of rank equal to $\sum_{i\ge
0}{\rm{dim}}_{\Q}H^{i}(X^{T};\Q)$.
\end{theorem}

\medskip

Recall that if $G$ acts on a topological space $X$ then its inertia
space is defined as $\Lambda X:=\{(g,x)\in G\times X ~|~ gx=x \}$.
Note that $\Lambda X$ is naturally a $G$--space via 
$h \cdot (x,g):=(hx, hgh^{-1})$. 
These spaces play an important role in geometry and
topology through the theory of stacks. A basic observation noted in
Section \ref{section max rank} is that if $G$ acts on $X$ with
connected maximal rank isotropy subgroups such that the
$\pi_1(G_{x})$ are torsion--free for all $x\in X$, then $\Lambda X$
also has connected maximal rank isotropy subgroups. In addition if
$X^T$ has the homotopy type of a finite $W$--CW complex the same is
true of $(\Lambda X)^T=X^T\times T$ and so we have:

\begin{theorem}
Let $X$ denote a compact $G$--CW complex with connected maximal rank
isotropy subgroups all of which have torsion--free fundamental
group. Then $K_G^*(\Lambda X)\otimes \Q$ is a free $R(G)\otimes
\Q$--module of rank equal to $2^{r}\cdot \left(\sum_{i\ge
0}{\rm{dim}}_{\Q}H^{i}(X^T;\Q)\right)$, where $r$ is the rank of $G$
as a compact Lie group.
\end{theorem}

The construction of inertia spaces can be iterated. In this way we
obtain a sequence of spaces $\{\Lambda^{n}(X)\}_{n\ge 0}$. Let $\P$
denote the collection of all compact Lie groups arising as finite
products of the classical groups $SU(r)$, $U(q)$ and $Sp(k)$. If we
further require that $G_{x}\in \P$ for all $x\in X$, then the action
of $G$ on $\Lambda^{n}(X)$ has connected maximal rank isotropy
subgroups for every $n\ge 0$ (Proposition \ref{action on inertia})
and given that $(\Lambda^n(X))^T= X^T\times T^n$, we have

\begin{theorem}\label{case of inertia spaces introduction}
Let $X$ denote a compact $G$--CW complex such that all of its
isotropy subgroups lie in $\P$ and are of maximal rank. Then
$K_G^*(\Lambda^n(X))\otimes \Q$ is a free $R(G)\otimes \Q$--module
of rank equal to $2^{nr}\cdot \left(\sum_{i\ge
0}{\rm{dim}}_{\Q}H^{i}(X^T;\Q)\right)$ where $r$ is the rank of $G$
as a compact Lie group.
\end{theorem}

Taking $X$ to be a single point with the trivial $G$--action yields
$\Lambda^n(X) = \Hom(\Z^n, G)$  (the space of ordered commuting $n$-tuples
in $G$) with the conjugation action, and our result can be applied.

\begin{corollary}\label{commuting elements introduction}
Suppose that $G\in \P$ is of rank $r$. Then
$K_{G}^{*}(\Hom(\Z^{n},G))\otimes \Q$ is free of rank $2^{nr}$ as an
$R(G)\otimes\Q$--module.
\end{corollary}

Now let $C_{n}(\g)^{+}$ denote the one point compactification of
the algebraic variety of ordered commuting $n$-tuples
in $\g$, the Lie algebra of $G$. This variety is endowed 
with an action of $G$ via the adjoint representation. 
As before we can show that if $G\in\P$ then
it has connected maximal rank isotropy and we obtain

\begin{corollary}\label{commuting variety intro}
Suppose that $G\in \P$ is of rank $r$. Then there is an isomorphism
of modules over $R(G)\otimes \Q$
\begin{equation*}
\tilde{K}^{q}_{G}(C_{n}(\g)^{+})\otimes \Q \cong \left\{
\begin{array}{ccc}
R(G)\otimes \Q &\text{ if } &q\equiv rn\ (\text{mod } 2),\\
0& \text{ if } &q+1\equiv rn\ (\text{mod } 2).
\end{array}%
\right.
\end{equation*}
\end{corollary}

Our methods also yield calculations for $K^*_G(X)$ but they are
somewhat more technical. To state them we need to recall a
definition from the theory of reflection groups. Let $\Phi$ denote
the root system associated to a fixed maximal torus $T\subset G$ and
let $\Phi^{+}$ be a choice of positive roots. Suppose that
$W_{i}\subset W$ is a reflection subgroup. Let $\Phi_{i}$ be the
corresponding root system and $\Phi_{i}^{+}$ the corresponding
positive roots. Define $W_{i}^{\ell}:=\{w\in W ~|~
w(\Phi^{+}_{i})\subset \Phi^{+}\}$. The set $W_{i}^{\ell}$ forms a
system of representatives for the different cosets in $W/W_{i}$. Let
$\W=\{W_{i}\}_{i\in \I}$ be a family of reflection subgroups of $W$.
We say that $\W$ satisfies the \textsl{coset intersection property}
if given $i,j\in \I$ we can find some $k\in \I$ such that $W_{i}\cup
W_{j}\subset W_{k}$ and $W_{k}^{\ell}=W_{i}^{\ell}\cap
W_{j}^{\ell}$. We are now ready to state our result:

\begin{theorem}\label{mainintroduction2}
Let $G$ be a compact connected Lie group with torsion--free
fundamental group and $T\subset G$ a maximal torus. Suppose that $X$
is a compact $G$--CW complex with connected maximal rank isotropy.
Assume that there is a CW-subcomplex $K$ of $X^{T}$ such that for
every element $x\in X^{T}$ there is a unique $w\in W$ such that
$wx\in K$ and the family $\{W_{\sigma}~~|~~ \sigma \text{ is a cell
in } K\}$ is contained in a family $\W$ of reflection subgroups of
$W$ satisfying the coset intersection property. If $H^{*}(X^{T};\Z)$
is torsion--free, then $K_{G}^{*}(X)$ is a free module over $R(G)$
of rank equal to $\sum_{i\ge 0}{\rm{rank}}_{\Z}H^{i}(X^{T};\Z)$.
\end{theorem}

\begin{corollary}\label{linear representations intro}
Let $G$ be a compact connected Lie group with $\pi_{1}(G)$
torsion--free. Let $G$ act on its Lie algebra $\g$ by the adjoint
representation.  If $r$ is the rank of $G$, then
\begin{equation*}
\tilde{K}^{q}_{G}(\S^{\g})\cong \left\{
\begin{array}{ccc}
R(G)&\text{ if } &q\equiv r\ (\text{mod } 2),\\
0& \text{ if } &q+1\equiv r\ (\text{mod } 2).
\end{array}%
\right.
\end{equation*}
\end{corollary}

\medskip

To prove our main results we show that under the appropriate
conditions the equivariant $K$--theory of a space with connected
maximal rank isotropy can be computed from the strong collapse of
the well--known spectral sequence associated to its skeletal
filtration which can be described using Bredon cohomology. This in
turn can be expressed in terms of the cohomology of the fixed point
set $X^T$; it requires an explicit identification of the $E_2$--term
of the spectral sequence as a module over the representation ring.
This involves some subtleties concerning the representation rings of
compact Lie groups related to work of R. Steinberg (see
\cite{Steinberg}).

\medskip

In \cite{Brylinski} it was established that if $G$ is compact and
connected with torsion--free fundamental group, then $K_G^*(G)$
(with the conjugation action) is a free $R(G)$--module of rank
$2^{r}$, where $r$ is the rank of $G$ (note that our methods provide
a different proof of this for $G$ simply connected, see Corollary
\ref{Bryl}). However the example of commuting pairs in $SU(2)$ which
is computed in Example \ref{example commuting tuples in SU(2)} shows
that the analogous result does not extend to spaces of commuting
elements without inverting the primes dividing the order of the Weyl
group $W$. Moreover the condition $G\in\P$ is required to ensure
connected maximal rank isotropy; for example this already fails for
$\Hom(\Z^2, Spin(7))$ even though $Spin(7)$ is a simply connected,
simple Lie group. This is due to the fact that $\Hom(\Z^{3},Spin(7))$ 
is not path--connected (see Example \ref{G on commuting} and 
\cite{KS}).

\medskip

The layout of this article is as follows. In Section \ref{section
max rank} we derive some basic properties of actions of compact Lie
groups that have connected maximal rank isotropy subgroups and
provide some key examples where such actions appear naturally. In
Section \ref{The representation ring} we construct a suitable basis
for $R(T)^{W_{i}}$, the ring of $W_{i}$ invariants in $R(T)$, for a
reflection subgroup $W_{i}$ of the Weyl group $W$. In Section
\ref{section Bredon cohomology} we provide a method for the
computation of the Bredon cohomology whose coefficient system is
obtained by taking fixed points of $R(T)$. In Section
\ref{equivariant K-theory} we prove Theorems \ref{mainintroduction1}
and \ref{mainintroduction2} which are the main results in this
article. In Section \ref{applications} some applications of the main
theorems are discussed. Finally the Appendix contains basic
definitions and facts about $G$--CW complexes and Bredon cohomology
which the reader may find useful.

\section{Group actions with maximal rank isotropy subgroups}
\label{section max rank}

Throughout this section $G$ will be a compact connected Lie group,
$T\subset G$ a maximal torus and $W$ the associated Weyl group.

\begin{definition}
Let $X$ be a $G$--space. The action of $G$ on $X$ is said to have 
maximal rank isotropy subgroups if for every $x\in X$, the isotropy 
group $G_{x}$ is a subgroup of maximal rank; that is, for every
$x\in X$ we can find a maximal torus $T_{x}$ in  $G$ such that
$T_{x}\subset G_{x}$. If in addition $G_{x}$ is connected for every
$x\in X$, then the action of $G$ on $X$ is said to have connected
maximal rank isotropy subgroups.
\end{definition}

Let $G$ act on a topological space $X$ with connected maximal
rank isotropy subgroups. Then by restriction the group $N_{G}(T)$ acts on
$X^{T}$ inducing an action of $W=WG:=N_{G}(T)/T$ on $X^{T}$. In the case
that $G$ acts smoothly on a smooth manifold $X$ the action of $G$ on $X$
is determined, up to equivalence, by the associated
action of $W$ on $X^{T}$ as was proved in \cite[Theorem 4.1]{Hauschild}.
In the broader context of \textsl{continuous} actions of a compact Lie
group on a topological space $X$ with connected maximal rank isotropy
subgroups, the corresponding action of $W$ on $X^{T}$ still determines
the action of $G$ on $X$ in some way. Given a subgroup $H$ of $G$
with $T\subset H$ let $WH=N_{H}(T)/T$ denote its Weyl group.
Let $\Fi_{G}(X)$ denote the set of conjugacy classes of isotropy subgroups
of the action of $G$ on $X$. The set $\Fi_{G}(X)$ is partially ordered by
inclusion. Similarly, let $\Fi_{W}(X^{T})$ be the partially ordered
set of conjugacy classes of isotropy subgroups of the action of $W$ on
$X^{T}$.

\medskip

Now let $X$, $Y$ be two spaces on
which $G$ acts with connected maximal rank isotropy subgroups.
Then the assignment
$\Fi_{G}(X)\to \Fi_{W}(X^{T})$ which sends
$(H)$ to $(WH)$
defines an order preserving bijection and
any $W$-equivariant continuous map $f:X^{T}\to Y^{T}$ has a unique
$G$-equivariant continuous extension $F:X\to Y$
(see \cite[Theorem 1.1]{Hauschild} and \cite[Theorem 2.1]{Hauschild}).

\medskip

The following theorem shows
that $G$--CW complex structures on a $G$--space
$X$ with connected maximal rank isotropy subgroups are determined,
up to equivalence, by $W$--CW complex structures on $X^{T}$.

\begin{theorem}\label{CW structures}
Suppose that $G$ acts on a compact space $X$ with connected maximal
rank isotropy subgroups. Then, up to equivalence, there is a one to
one correspondence between $G$--CW complex structures on $X$ with
isotropy type $\Fi_{G}(X)=\{(H)~|~ H\in \text{Iso}(X)\}$ and
$W$--CW complex structures on $X^{T}$ with isotropy type
$\Fi_{W}(X^{T})=\{(WH)~|~ H\in \text{Iso}(X)\}$.
\end{theorem}
\Proof Assume that $X$ has the structure of a $G$--CW complex. By
\cite[Theorem 2 1.14]{Dieck1} it follows that $X^{T}$ has a natural
structure as $W$--CW complex. Explicitly, suppose that
$\{X_{n}\}_{n\ge 0}$ is the skeletal filtration of $X$. Since $G$
acts on $X$ with connected maximal rank isotropy subgroups, it
follows that the cells in the $G$--CW decomposition of $X$ are
of the form $G/H\times \D^{n+1}$, where $H\subset G$ is connected
subgroup of maximal rank. The skeletal filtration of $X^{T}$ is
$\{(X^{n})^{T}\}_{n\ge 0}$ and $(X^{n+1})^{T}$ is obtained from
$(X^{n})^{T}$ by attaching $W/WH\times \D^{n+1}$ along the attaching
map $W/WH\times \S^{n}\to (X^{n})^{T}$ obtained by taking the
$T$-fixed points to the corresponding attaching $G$-map
$G/H\times \S^{n}\to X^{n}$. In particular, $X^{T}$ has orbit type
$\Fi_{W}(X^{T})=\{(WH)~|~ H\in \text{Iso}(X)\}$.

\medskip

Conversely, assume that $X^{T}$ has a structure of a $W$--CW
complex. Then $X^{T}$ is a compact $W$--CW complex, in particular,
$X^{T}$ is a $W$-ENR and by \cite[Remark 8.2.5]{Dieck} it satisfies
conditions ($W$--CW 1)-($W$--CW 3) (see the Appendix for the
definition). On the other hand, $(H)\mapsto (WH)$ defines an order
preserving bijection between $\Fi_{G}(X)$ and $\Fi_{W}(X^{T})$. This
shows that $\Fi_{W}(X^{T})=\{(WH)~|~ H\in \text{Iso}(X)\}$ and also
that $X$ has finite orbit type since $\Fi_{W}(X^{T})$ is finite.
Suppose that $H\in \text{Iso}(X)$. We may assume that $T\subset H$
after replacing $H$ with some conjugate if necessary. Let $(h,u)$ be
a representation of $((X^{T})^{(WH)},(X^{T})^{>(WH)})$ as a $W$-NDR
pair. Such a representation exists because $X^{T}$ satisfies
($W$--CW 2). Since the assignment $(H)\mapsto (WH)$ is an order
preserving bijection we have $(X^{T})^{(WH)}=(X^{(H)})^{T}$. Then
$h:(X^{(H)})^{T}\times I\to (X^{(H)})^{T}$ is a $W$-equivariant map
and the action of $G$ on $X^{(H)}$ has connected maximal rank
isotropy subgroups. Therefore there is a unique $G$-equivariant map
$\tilde{h}:X^{(H)}\times I \to X^{(H)}$ extending $h$. Similarly the
map $u$ can be extended to a $G$-invariant map $\tilde{u}:X^{(H)}\to
I$ and it is easy to see that the pair $(\tilde{h},\tilde{u})$ is a
representation of $(X^{(H)},X^{>(H)})$ as a $G$-NDR pair. On the
other hand, since $X^{T}$ satisfies ($W$--CW 3) we have
$X^{T}_{WH}\to X^{T}_{WH}/(N_{W}WH/WH)$ is a numerable principal
$N_{W}WH/WH$-bundle for every $H\in \text{Iso}(X)$. By \cite[Lemma
1.6]{Hauschild} we have $N_{G}(H)/H=N_{W}(WH)/WH$. Using that
$(H)\mapsto (WH)$ defines a bijection from $\Fi_{G}(X)$ to
$\Fi_{W}(X^{T})$ and \cite[Lemma 1.4]{Hauschild} we see that
$X_{H}=X^{T}_{WH}$ and thus $X_{H}\to X_{H}/N_{G}(H)/H$ is a
numerable principal $N_{G}(H)/H$-bundle. This proves that $X$
satisfies conditions ($G$--CW 1)-($G$--CW 3). By assumption $X^{T}$
is a $W$--CW complex, in particular Theorem \ref{Corollary Luck}
shows that $(X^{T})^{WH}$ has the homotopy type of a CW-complex for
all $H\in \text{Iso}(X)$. As before we may assume that $T\subset G$
after replacing $H$ with a suitable conjugate. In this case by
\cite[Lemma 1.4]{Hauschild} we have $X^{H}=(X^{T})^{WH}$ and this
space has the homotopy type of a CW-complex. Applying Theorem
\ref{Corollary Luck} again we see that the space $X$ can be given
the structure of a $G$--CW complex.

\medskip

Finally, an application of the equivariant Whitehead theorem shows
that, up to equivalence, this defines a one to one correspondence
between $G$--CW complex structures on $X$ and $W$--CW complex
structures on $X^{T}$.\qed

\medskip

Examples of spaces endowed with an action of $G$ with connected maximal
rank isotropy subgroups can be constructed in a natural way by
considering the action of $G$ on itself by conjugation as is
shown below.

\begin{example}\label{G on itself}
Suppose that $G$ is a compact connected Lie group with $\pi_{1}(G)$
torsion--free. Then the conjugation action of $G$ on itself has
connected maximal rank isotropy subgroups. To see this, note that
for every $g\in G$ the isotropy group under this action is
$G_{g}=Z_{G}(g)$, the centralizer of $g$ in $G$. When $\pi_{1}(G)$
is torsion--free $G_{g}$ is connected for every $g\in G$ by
\cite[Corollary IX 5.3.1]{Bourbaki}. Also, any element $g\in G$ is
contained in some maximal torus $T$, in particular $T\subset G_{g}$
and this means that $G_{g}$ is a maximal rank subgroup.
\end{example}

\begin{example}\label{G on commuting}
More generally let $n\ge 1$ be an integer and consider the space of
commuting $n$-tuples in $G$, $\Hom(\Z^{n},G)$. Let $\BONE:\Z^{n}\to G$ be
the trivial representation and denote the path--connected component of
$\Hom(\Z^{n},G)$ that contains $\BONE$ by $\Hom(\Z^{n},G)_{\BONE}$. This
component is precisely the subspace of $\Hom(\Z^{n},G)$ consisting of all
commuting $n$-tuples $(x_{1},\dots,x_{n})$ that are contained in some
maximal torus of $G$ by \cite[Lemma 4.2]{Baird}. The space
$\Hom(\Z^{n},G)$ can be seen as a $G$--space via the conjugation action of
$G$ and $\Hom(\Z^{n},G)_{\BONE}$ is invariant under this action. For any
$\ux=(x_{1},\dots,x_{n})\in \Hom(\Z^{n},G)$ the isotropy subgroup at $\ux$
under this action is $G_{\ux}=Z_{G}(\ux)$. By the previous comment, if
$\ux\in \Hom(\Z^{n},G)_{\BONE}$ then $Z_{G}(\ux)$ contains a maximal torus
of $G$. This shows that $\Hom(\Z^{n},G)_{\BONE}$ has maximal rank isotropy
subgroups. Moreover, the conjugation action of $G$ on $\Hom(\Z^{n},G)$ 
has connected maximal rank isotropy subgroups if and only if 
$\Hom(\Z^{n+1},G)$ is path--connected. Indeed, suppose first that 
$\Hom(\Z^{n+1},G)$ is path--connected. This implies in particular 
that  $\Hom(\Z^{n},G)_{\BONE}=\Hom(\Z^{n},G)$. Suppose that $z\in G$
commutes with $\ux=(x_{1},\dots,x_{n})\in \Hom(\Z^{n},G)_{\BONE}$. Then
\[
(z,x_{1},\dots,x_{n})\in \Hom(\Z^{n+1},G)=\Hom(\Z^{n+1},G)_{\BONE}
\]
and therefore  we can find a maximal torus $T$ of $G$ that contains
$z,x_{1},\dots,x_{n}$. In particular we can find a path $\gamma$ in $T$
joining $z$ and $1_{G}$. Note that $\gamma$ is contained in
$T\subset Z_{G}(\ux)$ proving that $Z_{G}(\ux)$ is path--connected. 
Conversely, suppose that the action of $G$ on $\Hom(\Z^{n},G)$ 
has connected isotropy subgroups. Assume that 
$(x_{1},\dots,x_{n+1})\in \Hom(\Z^{n+1},G)$. Then 
$x_{n+1}\in Z_{G}(x_{1},\dots,x_{n})$ and by assumption this is 
path--connected.
We conclude that we can find a path $\gamma$ in 
$Z_{G}(x_{1},\dots,x_{n})$ from $x_{n+1}$ to $1_{G}$. The path $\gamma$ 
shows that  $(x_{1},\dots,x_{n+1})$ lies in the same path--connected 
component of $(x_{1},\dots,x_{n},1_{G})$. Iterating this argument we 
conclude that $(x_{1},\dots,x_{n+1})$ lies in the same path--connected 
component of $\BONE=(1_{G},\dots,1_{G})$ proving that $\Hom(\Z^{n+1},G)$ 
is path--connected.

\end{example}

\medskip

The following proposition characterizes the compact Lie groups for which
$\Hom(\Z^{n},G)$ is path--connected for all $n\ge 0$ when $G$ is simple.

\begin{proposition}\label{characterization of connectedness}
If $G$ is a simple, simply connected compact Lie group, then the
following statements are equivalent:
\begin{enumerate}
\item $\Hom(\Z^n, G)$ is path--connected for all $n\ge 0$.
\item $A\subset G$ is a maximal abelian subgroup of $G$ if and only if $A$
is a maximal torus.
\item $G$ is isomorphic to either $SU(r)$ for some $r\ge 2$ or to $Sp(k)$
for some $k\ge 2$.
\end{enumerate}
\end{proposition}
\Proof We first prove that $(1)$ is equivalent to $(2)$. For this
equivalence we only need to assume that $G$ is a compact Lie group.
Suppose first that $\Hom(\Z^{n},G)$ is path--connected for all $n\ge
0$. Let $A$ be a maximal abelian subgroup of $G$. Note that $A$ must
be closed since $\overline{A}$ is also abelian. Therefore $A$ is a
compact Lie group. Let $A_{0}$ be the path--connected component of
$A$ containing the identity element. Then $A_{0}$ is a compact
connected abelian Lie group and thus it is a torus. In particular,
we can find a finite set $\{a_{1},\dots, a_{m}\}$ of elements in $A$
such that $A_{0}=\overline{\left<a_{1},\dots,a_{m}\right>}$. As $A$
has only finitely many distinct components, we may further choose
elements $b_{1},\dots b_{k}$ in $A-A_{0}$ such that
$A=\overline{\left<a_{1},\dots,a_{m},b_{1},\dots,b_{k}\right>}$.
Since $A$ is abelian we have
$(a_{1},\dots,a_{m},b_{1},\dots,b_{k})\in \Hom(\Z^{m+k},G)$. The
connectedness of the latter implies that there is a maximal torus
$T$ in $G$ such that $a_{i}, b_{j}\in T$ for all $i,j$. In
particular
$A=\overline{\left<a_{1},\dots,a_{m},b_{1},\dots,b_{k}\right>}\subset
T$. Since $A$ is maximal abelian this shows that $A=T$. Conversely,
suppose that $T$ is a maximal torus. If $T$ is not a maximal abelian
subgroup, then we can find some $g\in G$ with $g\notin T$ and such
that $g$ commutes with all the elements in $T$. Choose $a_{1},\dots,
a_{m}\in T$ such that $T=\overline{\left<a_{1},\dots,a_{m}\right>}$.
Note that $(a_{1},\dots,a_{m},g)\in \Hom(\Z^{m+1},G)$. Again by
connectedness of $\Hom(\Z^{m+1},G)$ there exists a maximal torus
$T'$ in $G$ such that $a_{1}\dots, a_{m}, g \in T'$. This implies
that $T=\overline{\left<a_{1},\dots,a_{m}\right>}\subsetneq T'$
which is a contradiction since $T$ is a maximal torus.

\medskip

Assume now that $(2)$ is true and suppose that $(g_{1},\dots g_{n})
\in \Hom(\Z^{n},G)$. Let $B=\left<g_{1},\dots g_{n}\right>$. Then $B$
is an abelian subgroup and we can find a maximal abelian group
$A\subset G$ with $B\subset A$. By assumption $A$ must be a torus,
in particular, $A$ is path--connected. It follows by
\cite[Proposition 2.3]{AC} that $\Hom(\Z^{n},G)$ is path--connected
for every $n\ge 0$. This proves that $(1)$ and $(2)$ are equivalent.

\medskip

By \cite[Corollary 2.4]{AC} if $G=SU(r)$ or $G=Sp(k)$ then
$\Hom(\Z^{n},G)$ is path--connected for every $n\ge 0$. This shows
that statement $(3)$ implies $(1)$. Assume now that $G$ is a simple,
simply connected compact Lie group. Table 2 in \cite{KS} together
with \cite[Theorem 3]{KS} show that under these assumptions
$\Hom(\Z^{3},G)$ is path--connected only when $G=SU(r)$ for some
$r\ge 2$ or $G=Sp(k)$ for some $k\ge 2$. This proves in particular
that $(1)$ implies $(3).$
\qed

\medskip

Let $G$ denote a compact simply connected Lie group. Then we may
express it as a product $G \cong G_1\times \dots \times G_s$, where
$G_1,\dots, G_s$ are simply connected simple Lie groups. Now it is
easy to see that $\Hom(\Z^n, H\times Q)\cong \Hom(\Z^n,H)\times
\Hom(\Z^n, Q)$ and hence $\Hom (\Z^n, H\times Q)$ is path--connected
if and only if $\Hom (\Z^n, H)$ and $\Hom (\Z^n, Q)$ are both
path--connected.

\begin{corollary}\label{cor connected hom}
Let $G$ be a compact simply connected Lie group, then
$\Hom(\Z^n, G)$ is path--connected for all $n\ge 1$ if
and only if $G$ is isomorphic to a finite cartesian product of
groups of the form $SU(r)$ ($r\ge 2$) and $Sp(k)$
($k\ge 2$).
\end{corollary}

We observe that the unitary groups $U(q)$ also satisfy the first two
conditions in the proposition above. This leads us to introduce the
following.

\begin{definition}
Let $\P$  denote the collection of all compact Lie groups arising as
finite cartesian products of the groups $SU(r)$, $U(q)$ and $Sp(k)$.
\end{definition}

Note that for $G\in \P$, $\Hom(\Z^{n},G)$ is path--connected for
every $n\ge 0$ and therefore as observed in Example \ref{G on
commuting} this $G$--space has connected maximal rank isotropy. We
now recall the notion of the inertia space for a group action.

\begin{definition}
Let $G$ denote a Lie group acting on a space $X$. Then its
inertia space is defined as
$\Lambda X = \{ (x, g) ~~|~~ gx=x\} \subset X\times G$.
\end{definition}

This construction can of course be iterated. This way we obtain a
sequence of spaces $\{ \Lambda^n(X)\}_{n\ge 1}$. Note that the
correspondence $(x,g)\mapsto (hx, hgh^{-1})$ defines a $G$--action
on $\Lambda X$, and the isotropy subgroups are simply $G_{(x,g)} =
Z_{G_x}(g)$, the centralizer of $g$ in $G_x$. Note that if $G$ acts
on $X$ with connected maximal rank isotropy subgroups with
torsion--free fundamental group, then the $G$--action on $\Lambda X$
will also have connected maximal rank isotropy by Example \ref{G on
itself}. Similarly we obtain that $G_{(x,g_1,\dots , g_n)} =
Z_{G_x}(g_1,\dots , g_n)$ where $(g_1, \dots , g_n)\in \Hom (\Z^n,
G_x)$. Note that $\Lambda^n(\{x_0\})=\Hom(\Z^n,G)$.

\begin{proposition}\label{action on inertia}
Let $G$ be a compact Lie group. Suppose $G$ acts on a space $X$ with
connected maximal rank isotropy subgroups in such a way that
$G_{x}\in \P$ for every $x\in X$. Then the induced action of $G$ on
$\Lambda^n(X)$ has connected maximal rank isotropy subgroups.
\end{proposition}
\Proof If $(x, g_1, \dots , g_n)\in \Lambda^n(X)$ then its isotropy
subgroup is $Z_{G_x}(g_1, \dots , g_n)$, which we have seen are
connected and of maximal rank. \qed

\medskip

Thus this approach provides a method for constructing natural spaces
on which $G$ acts with connected maximal rank isotropy subgroups.

\section{The representation ring}\label{The representation ring}

Suppose that $G$ is a compact connected Lie group and let $T\subset
G$ be a maximal torus. The inclusion map $i:T\to G$ induces a
natural isomorphism of rings
$i^{*}:R(G)\stackrel{\cong}{\rightarrow}R(T)^{W}$, see for example
\cite[IX 3]{Bourbaki}. We can see $R(T)$ as a module over $R(G)\cong
R(T)^{W}$. By \cite[Theorem 1]{Pittie} $R(T)$ is a free module over
$R(T)^{W}$ of rank $|W|$ when $\pi_{1}(G)$ is torsion--free. More
generally, by \cite[Theorem 2.2]{Steinberg} given a reflection
subgroup $W_{i}\subset W$ the ring of invariants $R(T)^{W_{i}}$ is a
free module over $R(T)^{W}$ of rank $|W|/|W_{i}|$. This is described
next; we follow the presentation provided in \cite{Uma} with some
modifications to fit our settings. Assume first that $G$ is simply
connected and thus semisimple. Let $\Lambda=X^{*}(T)=\Hom(T,\S^{1})$
be the character group of $T$. $\Lambda$ is a finitely generated
free abelian group (where the group operation is written
multiplicatively) and the representation ring $R(T)$ is the group
algebra $\Z[\Lambda]$. Let $\Phi$ be the root system associated to
$(G,T)$. Fix a subset $\Phi^{+}$ of positive roots of $\Phi$ and let
$\Delta=\{\alpha_{1},\dots,\alpha_{r}\}$ be an ordering of the
corresponding set of simple roots. Consider the fundamental weights
$\{\omega_{1},\dots \omega_{r}\}$ corresponding to the simple roots
$\{\alpha_{1},\dots,\alpha_{r}\}$. Since $G$ is assumed to be simply
connected the set $\{\omega_{1},\dots \omega_{r}\}$ forms a free
basis of $\Lambda$. Suppose that $W_{i}\subset W$ is a reflection
subgroup. Let $\Phi_{i}$ be the corresponding root system and
$\Phi_{i}^{+}$ the corresponding positive roots. Define
$W_{i}^{\ell}:=\{w\in W ~|~  w(\Phi^{+}_{i})\subset \Phi^{+}\}$. The
set $W_{i}^{\ell}$ forms a system of representatives of the
\textsl{left} cosets in $W/W_{i}$ by \cite[Lemma 2.5]{Steinberg}.
This means that any element $w\in W$ can be factored in a unique way
in the form $w= ux \text{ with } u\in W_{i}^{\ell} \text{ and } x\in
W_{i}$. Moreover, when $\Phi_{i}$ has a basis consisting of all
roots of the form $\{\alpha~|~ \alpha\in I\}$ for some subset
$I\subset \Delta$, then $W_{i}^{\ell}=W^{I}:=\{w\in W ~|~
\ell(ws)>\ell(w) \text{ for all } s\in I\}$ which is the set of
minimal length coset representatives of $W/W_{i}$ with respect to
the length function $\ell$ determined by $\Delta$. Because of this,
in general we call the set $W_{i}^{\ell}$ a system of minimal length
representatives of $W/W_{i}$. Given $v\in W_{i}^{\ell}$ define
$p_{v}:=\prod_{v^{-1}\alpha_{i}<0}\omega_{i}$ and
$f^{W_{i}}_{v}:=\sum_{x\in W_{i}(v)\setminus
W_{i}}x^{-1}v^{-1}p_{v}$. Here $W_{i}(v):=\{w\in W_{i} ~|~
w^{-1}v^{-1}p_{v}=v^{-1}p_{v}\}$ is the stabilizer in $W_{i}$ of
$v^{-1}p_{v}$. Thus $f^{W_{i}}_{v}$ is the sum of the elements in
the orbit of $v^{-1}p_{v}$ under the action of $W_{i}$, making
$f^{W_{i}}_{v}$ invariant under the action of $W_{i}$. By
\cite[Theorem 2.2]{Steinberg} the collection
$\{f^{W_{i}}_{v}\}_{v\in W_{i}^{\ell}}$ forms a free basis of
$R(T)^{W_{i}}$ as a module over $R(T)^{W}$. We refer to this basis
as the Steinberg basis. Next we discuss the relationship between the
Steinberg bases associated to reflection subgroups $W\supset
W_{j}\supset W_{i}$. We set up some notation first.

\medskip

Given $v\in W_{j}^{\ell}$, the subgroup $W_{j}(v)$ is itself a
reflection subgroup. Let $\Phi_{j}(v)$ be its corresponding root
system and $\Phi_{j}^{+}(v)$ the corresponding positive roots.
Define $W_{j}^{r}(v):=\{w\in W_{j} ~|~
w^{-1}(\Phi_{j}^{+}(v))\subset \Phi_{j}^{+}\}$. As before, this set
gives a system of representatives of the \textsl{right} cosets in
$W_{j}(v)\setminus W_{j}$ and it has the property that every element
$w\in W_{j}$ can be factorized in a unique way in the form $w= xu
\text{ with } x\in W_{j}(v) \text{ and } u\in W_{j}^{r}(v)$. Note
that by definition of $f_{v}^{W_{j}}$ we have
$f^{W_{j}}_{v}=\sum_{x\in W_{j}(v)\setminus W_{j}}x^{-1}v^{-1}p_{v}
=\sum_{x\in W_{j}^{r}(v)}x^{-1}v^{-1}p_{v}$. On the other hand,
given $v\in W^{\ell}_{j}$ and $x\in W_{j}^{r}(v)$ we can consider
the coset $xW_{i}\in W_{j}/W_{i}$. Let $m_{j,i}(x)$ be the minimal
length representative of the coset $xW_{i}\in W_{j}/W_{i}$.
Precisely, $m_{j,i}(x)$ is the unique element in $[W_{i}^{\ell}\
]^{W_{j}}:=\{w\in W_{j} ~|~  w(\Phi_{i}^{+})\subset \Phi_{j}^{+}\}$
such that $xW_{i}=m_{j,i}(x)W_{i}$. The following lemma describes
the precise relationship between the bases $\{f^{W_{j}}_{v}\}_{v\in
W_{j}^{\ell}}$ and $\{f^{W_{i}}_{v}\}_{v\in W_{i}^{\ell}}$.

\begin{lemma}\label{relationship}
Suppose that $W\supset W_{j}\supset W_{i}$ are reflection subgroups.
Then $W_{j}^{\ell}\subset W_{i}^{\ell}$ and if $v\in W_{j}^{\ell}$
then $f_{v}^{W_{j}} =\sum_{\tilde{x}\in
m_{j,i}(W_{j}^{r}(v))}f^{W_{i}}_{v\tilde{x}}$.
\end{lemma}
\Proof It is clear from the definition of $W_{j}^{\ell}$ that the
inclusion $W_{j}\supset W_{i}$ implies $W_{j}^{\ell}\subset
W_{i}^{\ell}$. Note that the right hand side of the equation in
Lemma \ref{relationship} is well defined because given $v\in
W_{j}^{\ell}$ and $\tilde{x}\in m_{j,i}(W_{j}^{r}(v))$ then
$v\tilde{x}\in W_{i}^{\ell}$. We are going to prove first the
assertion for the particular case $W_{i}=\{1\}$. In this case
$x=m_{j,i}(x)$ for every $x\in W_{j}^{r}(v)$, thus we need to prove
\begin{equation}\label{reduction}
f_{v}^{W_{j}}=\sum_{x\in W_{j}^{r}(v)}f^{\{1\}}_{vx}.
\end{equation}
Suppose that  $v\in W_{j}^{\ell}$ and $x\in W_{j}^{r}(v)$. An argument
similar to the one provided in \cite[Proposition 1.9]{Uma} shows that
under these hypotheses
\begin{equation}\label{identity}
x^{-1}v^{-1}p_{v}=x^{-1}v^{-1}p_{vx}.
\end{equation}
Then
\[
f_{v}^{W_{j}}=\sum_{x\in W_{j}^{r}(v)}x^{-1}v^{-1}p_{v}
=\sum_{x\in W_{j}^{r}(v)}x^{-1}v^{-1}p_{vx}
=\sum_{x\in W_{j}^{r}(v)}f^{\{1\}}_{vx}
\]
proving the lemma when $W_{i}=\{1\}$. We now show that the result is
true for any reflection subgroups $W\supset W_{j}\supset W_{i}$. Suppose
that $v\in W_{j}^{\ell}$. Then by (\ref{reduction})
\[
f_{v}^{W_{j}}=\sum_{x\in W_{j}^{r}(v)}f^{\{1\}}_{vx}
=\sum_{x\in W_{j}^{r}(v)}x^{-1}v^{-1}p_{vx}.
\]
For any $x \in W_{j}^{r}(v)$ we have a unique factorization
$x=m_{j,i}(x)y(x)$ with $m_{j,i}(x)\in [W_{i}^{\ell}\ ]^{W_{j}}$ and
$y(x)\in W_{i}$. Thus
\begin{equation}\label{sum1}
f_{v}^{W_{j}}=\sum_{x\in W_{j}^{r}(v)}x^{-1}v^{-1}p_{vx}
=\sum_{x\in W_{j}^{r}(v)}y(x)^{-1}m_{j,i}(x)^{-1}v^{-1}p_{vm_{j,i}(x)y(x)}.
\end{equation}
We can rearrange the terms in the previous sum in the following way.
Let $\tilde{x}=m_{j,i}(x)$, then $x=\tilde{x}y(x)$. When $\tilde{x}$
runs through all elements in $m_{j,i}(W_{j}^{r}(v))$ the set of values of
$y$ such that $y^{-1}\tilde{x}^{-1}v^{-1}p_{v\tilde{x}yx}$
is a term in the sum (\ref{sum1}) is precisely $W_{i}^{r}(v\tilde{x})$.
This shows that
\begin{align*}
f_{v}^{W_{j}}&=\sum_{x\in W_{j}^{r}(v)}y(x)^{-1}
m_{j,i}(x)^{-1}v^{-1}p_{vm_{j,i}(x)y(x)}\\
&=\sum_{\tilde{x}\in m_{j,i}(W_{j}^{r}(v))}
\left(\sum_{y\in W_{i}^{r}(v\tilde{x})} y^{-1}\tilde{x}^{-1}v^{-1}
p_{v\tilde{x}y}\right)\\
&=\sum_{\tilde{x}\in m_{j,i}(W_{j}^{r}(v))}
f^{W_{i}}_{v\tilde{x}}.
\end{align*}
The last equality is obtained by applying (\ref{reduction}).
\qed

\medskip

For our applications it will be more convenient to obtain a different
basis of $R(T)^{W_{i}}$ as a module over $R(T)^{W}$ that behaves better
under the natural inclusion $R(T)^{W_{j}}\subset R(T)^{W_{i}}$, whenever
$W\supset W_{j}\supset W_{i}$ belong to a suitable
family $\W=\{W_{i}\}_{i\in \I}$ of reflection subgroups of $W$.
The families of subgroups that we consider satisfy the following
condition.

\begin{definition}\label{coset intersection property}
Let $\W=\{W_{i}\}_{i\in \I}$ be a family of reflection subgroups of $W$.
We say that $\W$ satisfies the coset intersection property if given
$i,j\in \I$ we can find some $k\in \I$ such that
$W_{i}\cup W_{j}\subset W_{k}$ and
$W_{k}^{\ell}=W_{i}^{\ell}\cap W_{j}^{\ell}$.
\end{definition}

\begin{example}\label{example coset} Suppose we have a sequence
of reflection subgroups
\[
W=W_{0}\supset W_{1}\supset W_{2}\supset\cdots\supset W_{k}=\{1\}.
\]
Then by passing to the minimal coset representatives we obtain an
increasing sequence
\[
\{1\}=W_{0}^{\ell}\subset W^{\ell}_{1}\subset W^{\ell}_{2}\subset
\cdots\subset W^{\ell}_{k}=W.
\]
This shows that the family $\W=\{W_{i}\}_{0\le i\le k}$ satisfies
the coset intersection property. We are also going to consider the
following important example. As before let
$\Delta=\{\alpha_{1},\dots,\alpha_{r}\}$ be an ordering of the
corresponding set of simple roots associated to a set $\Phi^{+}$ of
positive roots of the root system $\Phi$. For every $I\subset
\Delta$ let $W_{I}$ be the reflection subgroup of $W$ generated by
the corresponding reflections $s_{\alpha}$ with $\alpha\in I$. The
family $\W=\{W_{I}\}_{I\subset \Delta}$ satisfies the coset
intersection property. To see this, suppose that $I,J\subset
\Delta$. Then $W_{I}\cup W_{J}\subset W_{I\cup J}$ and
$W^{\ell}_{I\cup J}=W_{I}^{\ell}\cap W_{J}^{\ell}$ because given any
$I\subset \Delta$ we have
$$W_{I}^{\ell}=W^{I}:=\{w\in W ~|~
\ell(ws)>\ell(w) \text{ for all } s\in I\}.$$
\end{example}

\bigskip

Now fix $\W=\{W_{i}\}_{i\in \I}$ a family of reflection subgroups
of $W$ satisfying the coset intersection property. Consider the family
$\W^{\ell}:=\{W^{\ell}_{i}\}_{i\in \I}$ of minimal length coset
representatives. For every $i\in \I$ define
\[
C_{i}=W^{\ell}_{i}\setminus
\left(\bigcup_{j\in \I, W_{i}\subsetneq W_{j}}W_{j}^{\ell} \right).
\]
The different sets of the form $C_{i}$  provide a decomposition of
$W$ into disjoint sets $W=\bigsqcup_{i\in \I}C_{i}$. This is
precisely where the coset intersection property is used. Indeed, any
element $v\in W$ is contained in some $C_{i}$ for some $i\in \I$ and
if $v\in C_{i}\cap C_{j}$, then $v\in W_{i}^{\ell}\cap W_{j}^{\ell}$
and by assumption we can find some $k\in \I$ such that $W_{i}\cup
W_{j} \subset W_{k}$ and $W_{k}^{\ell}=W_{i}^{\ell}\cap
W_{j}^{\ell}$. If $W_{i}\subsetneq W_{k}$ we conclude that $v\notin
C_{i}$ which contradicts our assumption. A similar situation occurs
if $W_{j}\subsetneq W_{k}$. Therefore $W_{k}=W_{i}=W_{k}$ and thus
$i=j=k$.

\medskip

Additionally, for any $i\in \I$ we have
\[
W_{i}^{\ell}=\bigsqcup_{\substack{j\in \I\\ W_{i}\subset W_{j}}}C_{j}
\]
Thus given any $v\in W$ we can find a unique $i(v)\in \I$ such that
$v\in C_{i(v)}^{\ell}$ and we can define $g_{v}:=f^{W_{i(v)}}_{v}$.

\begin{theorem}\label{free basis}
With the previous notation, if $G$ is simply connected and
$\W=\{W_{i}\}_{i\in \I}$ is a family of reflection subgroups of $W$
satisfying the coset intersection property, then for every $i\in I$ the
collection $\{g_{v}\}_{v\in W_{i}^{\ell}}$ is a free basis of
$R(T)^{W_{i}}$ as a module over $R(T)^{W}$.
\end{theorem}
\Proof Let $i\in \I$, we will show that $\{g_{v}\}_{v\in W_{i}^{\ell}}$
is a free basis of $R(T)^{W_{i}}$ as a module over $R(T)^{W}$. Since
$R(T)^{W}$ is a domain, it suffices to show that every element of
$\{g_{v}\}_{v\in W_{i}^{\ell}}$ is in the span of
$\{f^{W_{i}}_{v}\}_{v\in W_{i}^{\ell}}$ and vice versa. We will prove
this by induction on
\[
n(i):=|\{j\in \I~~|~~ W_{i}\subset W_{j}\}|.
\]
When $n(i)=1$ we have $W_{i}=W$ and thus $W_{i}^{\ell}=\{1\}$. In this
case the claim is trivial as $g_{1}=f^{W}_{1}=1$. Suppose that the claim
is true for all $i\in \I$ for which $n(i)<n$. Let $i\in I$ be such that
$n(i)=n$. Suppose that $v\in W_{i}^{\ell}=
\bigsqcup_{\substack{j\in \I\\ W_{i}\subset W_{j}}}C_{j}$.
We show first that $g_{v}$ is in the span of
$\{f^{W_{i}}_{w}\}_{w\in W_{i}^{\ell}}$. If $v\in C_{i}$ then there is
nothing to prove as $g_{v}=f_{v}^{W_{i}}$. Suppose then that
$v\in C_{j}$ for some $j\in \I$ such that $W_{i}\subsetneq W_{j}$. Then
$v\in W_{j}^{\ell}\subset W_{i}^{\ell}$ and by Lemma \ref{relationship}
we have
\[
g_{v}=f^{W_{j}}_{v}
=\sum_{\tilde{x}\in m_{j,i}(W_{j}^{r}(v))}f^{W_{i}}_{v\tilde{x}}
\]
which is what we needed to show. Conversely, let's show that if
$v\in W_{i}^{\ell}$ then $f^{W_{i}}_{v}$ is in the span of
$\{g_{v}\}_{v\in W_{i}^{\ell}}$. If $v\in C_{i}^{\ell}$ then
$f_{v}^{W_{i}}=g_{v}$ and there is nothing to prove. Suppose then that
$v\in C_{j}$ for some $j\in \I$ such that $W_{i}\subsetneq W_{j}$.
Using Lemma \ref{relationship} we have
\begin{equation}\label{equax1}
f_{v}^{W_{j}}=\sum_{\tilde{x}\in m_{j,i}(W_{j}^{r}(v))}
f^{W_{i}}_{v\tilde{x}}
=f_{v}^{W_{i}}+\sum_{\substack{\tilde{x}\in m_{j,i}(W_{j}^{r}(v))\\
\tilde{x}\ne 1}}
f^{W_{i}}_{v\tilde{x}}.
\end{equation}
In this case $g_{v}=f_{v}^{W_{j}}$. Also, the inductive hypothesis
shows that for all values of $\tilde{x}$ such that $v\tilde{x}\in C_{k}$
for some $k\in \I$ with $W_{i}\subsetneq W_{k}$ then
$f^{W_{i}}_{v\tilde{x}}$ is in the span of
$\{g_{v}\}_{v\in W_{k}^{\ell}}\subset \{g_{v}\}_{v\in W_{i}^{\ell}}$.
Therefore (\ref{equax1}) can be rewritten in the form
\[
f^{W_{i}}_{v}=\sum_{v_{1}\in I_{1}}r_{v_{1}}g_{v_{1}}-
\sum_{\substack{\tilde{x}_{1}\in m_{j,i}(W_{j}^{r}(v))\\
\tilde{x}_{1}\ne 1, v\tilde{x}_{1}\in C_{i}}}
f^{W_{i}}_{v\tilde{x}_{1}}
\]
for some set $I_{1}\subset W_{i}^{\ell}$ and some $r_{v_{1}}\in
R(T)^{W}$. If there are no elements with $\tilde{x}_{1}\in
m_{j,i}(W_{j}^{r}(v))$ such that $v\tilde{x}_{1}\in C_{i}$ we are
done. Otherwise, we iterate this procedure. At each step we find a
sequence of elements $\tilde{x}_{i}\ne 1$ for $1\le i\le k$ such
that $\{v,v\tilde{x}_{1},\dots, v(\tilde{x}_{1}...\tilde{x}_{k})\}
\subset C_{i}$. Since $C_{i}$ is a finite set this process must
terminate after finitely many steps and thus $f_{v}^{W_{i}}$ is in
the span of $\{g_{v}\}_{v\in W_{i}^{\ell}}$. \qed

\bigskip

The important additional property of the basis constructed in the
previous theorem is that whenever $W_{j}\supset W_{i}$ then
$W_{j}^{\ell}\subset W_{i}^{\ell}$ and the set
$\{g_{v}\}_{v\in W_{j}^{\ell}}$ is contained in the set
$\{g_{v}\}_{v\in W_{i}^{\ell}}$. With this choice of basis we can
find an isomorphism of $R(T)^{W}$--modules for every $i\in \I$
\[
\varphi_{i}:R(T)^{W_{i}}\stackrel{\cong}
{\rightarrow} \bigoplus_{v\in W_{i}^{\ell}}
R(T)^{W}g_{v}
\]
that satisfies the following compatibility condition: whenever
$W_{j}\supset W_{i}$ the following diagram commutes
\begin{equation*}
\begin{CD}
R(T)^{W_{j}}@>>>R(T)^{W_{i}}\\
@V\varphi_{j}VV     @VV\varphi_{i}V\\
\bigoplus_{v\in W_{j}^{\ell}}
\displaystyle{R(T)^{W}g_{v}}@>>>
\bigoplus_{v\in W_{i}^{\ell}}
R(T)^{W}g_{v}.\\
\end{CD}
\end{equation*}
In this diagram the horizontal maps are the natural inclusions.

\medskip

Suppose now that $G$ is a compact connected Lie group with
$\pi_{1}(G)$ torsion--free. Then as pointed out in \cite{Pittie}
there is a finite covering sequence
\begin{equation}\label{coverspace}
1\to \Gamma\to T'\times G_{0}\stackrel{\pi}{\rightarrow} G\to 1,
\end{equation}
with $G_{0}$ a simply connected compact Lie group, $T'$ a torus and
$\Gamma$ a finite central subgroup. Note that
$\pi^{-1}(T)=T^{'}\times T_{0}$, where $T_{0}\subset G_{0}$ is a
maximal torus and there is a covering space
\[
1\to \Gamma\to T'\times T_{0}\to T\to 1.
\]
Let $W_{0}$ be the Weyl group  associated to $(G_{0},T_{0})$. The
covering sequence (\ref{coverspace}) shows that $W_{0}=W$. Let
$p_{1}:T'\times G_{0}\to T'$ be the projection map and
$\bar{\Gamma}=p_{1}(\Gamma)\subset T'$. For every  $i\in \I$ there
is an isomorphism $R(T)^{W_{i}}\cong R(T_{0})^{W_{i}}\otimes
R(T'/\bar{\Gamma})$. In particular, if $\{g_{v}\}_{v\in
W_{i}^{\ell}}$ is the free basis of $R(T_{0})^{W_{i}}$ as a module
over $R(T_{0})^{W}$ constructed above, then the collection
$\{g_{v}\otimes 1\}_{v\in W_{i}^{\ell}}$ is a free basis of
$R(T)^{W_{i}}$ as a module over $R(T)^{W}$. By abuse of notation we
also denote this basis by $\{g_{v}\}_{v\in W_{i}^{\ell}}$. The
following theorem summarizes the above.

\begin{theorem}\label{Pittie}
Let $G$ be a compact, connected Lie group with $\pi_{1}(G)$ torsion--free.
Let $\W=\{W_{i}\}_{i\in \I}$ be a family of reflection subgroups of $W$
satisfying the coset intersection property. Then for every $i\in \I$
there is an isomorphism of $R(T)^{W}$--modules
\[
\varphi_{i}:R(T)^{W_{i}}\stackrel{\cong}
{\rightarrow} \bigoplus_{v\in W_{i}^{\ell}}
R(T)^{W}g_{v}.
\]
These isomorphisms satisfy the following compatibility condition:
whenever $W_{j}\supset W_{i}$ the following diagram commutes
\begin{equation*}
\begin{CD}
R(T)^{W_{j}}@>>>R(T)^{W_{i}}\\
@V\varphi_{j}VV     @VV\varphi_{i}V\\
\bigoplus_{v\in W_{j}^{\ell}}
\displaystyle{R(T)^{W}g_{v}}@>>>
\bigoplus_{v\in W_{i}^{\ell}}
R(T)^{W}g_{v}.\\
\end{CD}
\end{equation*}
In this diagram the horizontal maps are the natural inclusions.
\end{theorem}

\medskip

\begin{example} Suppose that $G=SU(3)$. Let $T\subset SU(3)$ be
the maximal torus consisting of those $3$ by $3$ diagonal matrices with
entries $x_{1},x_{2}$ and $x_{3}$ in $\S^{1}$ such that
$x_{1}x_{2}x_{3}=1$. The Weyl group $W=\Sigma_{3}$ acts by permutation
of the diagonal entries in $T$. We can see the entries of an element
in $T$ as linear characters $x_{i}:T\to \S^{1}$ that satisfy the
equation $x_{1}x_{2}x_{3}=1$ and
\[
R(T)=\Z[x_{1},x_{2},x_{3}]/(x_{1}x_{2}x_{3}=1).
\]
The roots associated to the pair $(G,T)$ are $x_{i}x_{j}^{-1}$ for
$i\ne j$.  In this case $\Delta=\{x_{1}x_{2}^{-1},x_{2}x_{3}^{-1}\}$ is
a set of simple roots. The fundamental weights associated to this system
of simple roots are $\{x_{1},x_{1}x_{2}\}$ respectively. Enumerate the
elements in $W$ as follows
\[
v_{1}=1, v_{2}=(23),v_{3}=(123),v_{4}=(132),v_{5}=(12), v_{6}=(13).
\]
Consider the family $\W=\{W_{0},W_{1},W_{2}\}$ of reflection subgroups
of $W$, where
\[
W=W_{0}\supset W_{1}=\left\langle (12)\right\rangle\supset W_{2}=\{1\}.
\]
As pointed out before this sequence satisfies the coset intersection
property and we have
\[
\{1\}=\{v_{1}\}=W^{\ell}_{0}\subset W_{1}^{\ell}=\{v_{1},v_{2},v_{3}\}
\subset W_{2}^{\ell}=W.
\]
The different Steinberg bases $\{f_{v}^{W_{i}}\}_{v\in
W^{\ell}_{i}}$ of $R(T)^{W_{i}}$ as a module over $R(T)^{W}$ are
given below:
\begin{itemize}
\item For $W_{0}=W$ we have the trivial basis $f_{v_{1}}^{W_{0}}=1$.

\item For $W_{1}=\left\langle (12)\right\rangle$ we obtain the basis
$f_{v_{1}}^{W_{1}}=1,\ f_{v_{2}}^{W_{1}}=x_{1}x_{3}+x_{2}x_{3},\
f_{v_{3}}^{W_{1}}=x_{3}$.

\item For $W_{2}=\{1\}$ the corresponding basis is
\[
f_{v_{1}}^{W_{2}}=1,\ f_{v_{2}}^{W_{2}}=x_{1}x_{3},\
f_{v_{3}}^{W_{2}}=x_{3},\ f_{v_{4}}^{W_{2}}=x_{2}x_{3},\
f_{v_{5}}^{W_{2}}=x_{2},\ f_{v_{6}}^{W_{2}}=x_{2}x_{3}^{2}.
\]
\end{itemize}
On the other hand, the different bases $\{g_{v}\}_{v\in
W^{\ell}_{i}}$ of $R(T)^{W_{i}}$ as a module over $R(T)^{W}$ are as
follows:
\begin{itemize}
\item For $W_{0}=W$ we have the trivial basis $g_{v_{1}}=1$.
\item For $W_{1}=\left\langle (12)\right\rangle$ we obtain the basis
$g_{v_{1}}=1,\ g_{v_{2}}=x_{1}x_{3}+x_{2}x_{3},\ g_{v_{3}}=x_{3}$.
\item For $W_{2}=\{1\}$ the corresponding basis is
\[
g_{v_{1}}=1,\ g_{v_{2}}=x_{1}x_{3}+x_{2}x_{3},\
g_{v_{3}}=x_{3},\ g_{v_{4}}=x_{2}x_{3},\
g_{v_{5}}=x_{2},\ g_{v_{6}}=x_{2}x_{3}^{2}.
\]
\end{itemize}
\end{example}

\section{Bredon cohomology}\label{section Bredon cohomology}

In this section we use the bases obtained in the previous section to
provide a computation of certain Bredon cohomology groups; a brief
description on this invariant can be found in the Appendix.

\medskip

To start assume that $G$ is a compact connected Lie group. Fix
$T\subset G$ a maximal torus and let $W$ be the corresponding Weyl
group. Consider the $W$--modules $R(T)$ and $R(G)\otimes \Z[W]$,
where $\Z[W]$ denotes the group ring. These modules induce coefficient
systems
\[
\lR_{T}:=H^{0}(-,R(T)) \text{ and  } \lZ_{W}:=H^{0}(-,R(G)\otimes \Z[W]).
\]
Explicitly the values of these coefficient systems at an orbit of
the form $W/W_{i}$ are
$$\lR_{T}(W/W_{i})=R(T)^{W_{i}} \text{ and }
\lZ_{W}(W/W_{i})=R(G)\otimes \Z[W]^{W_{i}}\cong R(G)\otimes
\Z[W/W_{i}].$$
For every subgroup $W_{i}\subset W$ the abelian
groups $R(T)^{W_{i}}$ and $R(G)\otimes \Z[W/W_{i}]$ have the
structure of a module over $R(G)\cong R(T)^{W}$. In particular,
$\lR_{T}$ and $\lZ_{W}$ are coefficient systems in the category of
$R(G)$--modules. Suppose that $Y$ is a $W$--CW complex. Note that in
particular $Y$ has the structure of a CW-complex by forgetting the
action of $W$.

\begin{theorem}\label{Bredon cohomology}
Suppose that $G$ is a compact connected Lie group with $\pi_{1}(G)$
torsion--free. Let $Y$ be a $W$--CW complex with cells of the form
$W/W_{i}\times \D^{n}$ for reflection subgroups $W_{i}\subset W$.
Assume that there is a CW-subcomplex $K$ of $Y$ such that  for every
$x\in Y$ there is a unique $w\in W$ such that $wx\in K$ and the
family $\{W_{\sigma}~~|~~ \sigma \text{ is a cell in }  K\}$ is
contained in a family $\W$ of reflection subgroups of $W$ satisfying
the coset intersection property. Then there are isomorphisms of
$R(G)$--modules $H^{*}_{W}(Y;\lR_{T})\cong H^{*}_{W}(Y;\lZ_{W})\cong
H^{*}(Y;R(G))$.
\end{theorem}
\Proof Consider the CW-complex structure on $Y$ obtained by
forgetting the action of $W$. Let $K$ be a sub CW-complex of $Y$
such that every element in $Y$ is conjugated to a unique element in
$K$ and let be a family $\W=\{W_{i}\}_{i\in \I}$ of reflection
subgroups of $W$ containing the subgroups $W_{\sigma}$ for all cells
$\sigma$ in $K$ and satisfying the coset intersection property. Then
we can use Theorem \ref{Pittie} to find a free basis
$\{g_{v}\}_{v\in W_{i}^{\ell}}$ of $R(T)^{W_{i}}$ as a module over
$R(T)^{W}\cong R(G)$ for every $i\in \I$ that is compatible with the
inclusions $W_{j}\supset W_{i}$ for $i,j\in \I$. We can construct in
a similar way a free basis of $R(G)\otimes \Z[W]^{W_{i}}$ as a
module over $R(G)$ for every $i\in \I$ in the following way.  Given
$v\in W_{i}^{\ell}$ define $l_{v}^{W_{i}}:=1\otimes(\sum_{x\in
W_{i}}x^{-1}v^{-1})\in R(G)\otimes \Z[W]^{W_{i}}$. The collection
$\{l_{v}^{W_{i}}\}_{v\in W^{\ell}_{i}}$ forms a free basis of
$R(G)\otimes \Z[W]^{W_{i}}$ as a module over $R(G)$. Using this
basis we can construct a new basis that behaves better under the
different inclusions $W_{j}\supset W_{i}$ with $i,j\in \I$ as
follows. Given $v\in W$ we can find a unique $i(v)\in \I$ such that
$v\in C_{i(v)}$. Define $m_{v}=l_{v}^{W_{i(v)}}$. An argument
similar to that of Theorem \ref{free basis} shows that
$\{m_{v}\}_{v\in W_{i}^{\ell}}$ is a free basis of $R(G)\otimes
\Z[W]^{W_{i}}$ as a module over $R(G)$. These bases enjoy the
further property that $\{m_{v}\}_{v\in W_{j}^{\ell}}$ is a subset of
$\{m_{v}\}_{v\in W_{i}^{\ell}}$ whenever $W_{j}\supset W_{i}$ for
$i,j\in \I$. Note that the bases $\{g_{v}\}_{v\in W_{i}^{\ell}}$ and
$\{m_{v}\}_{v\in W_{i}^{\ell}}$ provide an isomorphism of
$R(G)$--modules for every $i\in\I$
\begin{align*}
\psi_{W_{i}}:R(T)^{W_{i}}&\to R(G)\otimes \Z[W]^{W_{i}}\\
g_{v}&\mapsto m_{v}
\end{align*}
in such a way that whenever $W_{j}\supset W_{i}$ for $i,j\in \I$ the
following diagram commutes
\begin{equation}\label{commute iso}
\begin{CD}
R(T)^{W_{j}}@>>>R(T)^{W_{i}}\\
@V\psi_{W_{j}}VV     @VV\psi_{W_{i}}V\\
\bigoplus_{v\in W_{i}^{\ell}}
R(G)\otimes \Z[W]^{W_{j}}@>>>
\bigoplus_{v\in W_{j}^{\ell}}
R(G)\otimes \Z[W]^{W_{i}}.\\
\end{CD}
\end{equation}
In the previous diagram the horizontal maps are the inclusion maps;
we now show that the different isomorphisms $\{\psi_{W_{i}}\}_{i\in
\I}$ induce an isomorphism of $R(G)$--modules
$H^{*}_{W}(Y;\lR_{T})\cong H^{*}_{W}(Y;\lZ_{W})$. In fact, we are
going to show that the cochain complexes computing these Bredon
cohomology groups are isomorphic. To see this recall that
\[
C_{W}^{n}(Y;\lR_{T})=\bigoplus_{\sigma\in S_{n}(Y)} R(T)^{W_{\sigma}},
\]
where $S_{n}(Y)$ is a set of representatives of all $W$-cells in $Y$.
For $n\ge 0$ let $I_{n}(K)$ be the set of all $n$-dimensional cells in
$K$. Since every element in $Y$ is conjugated to a
unique element in $K$ it follows that $K$ has a unique representative
for all $W$-cells in $Y$,  thus we can choose $S_{n}(Y)$ to be
$I_{n}(K)$. As explained in the Appendix given $x\in C_{W}^{n}(Y;\lR_{T})$
and any $\sigma\in I_{n+1}(K)$
\[
\delta(x)_{\sigma}=\sum_{\tau\in I_{n}(K)}[\tau:\sigma]
i^{*}_{\tau,\sigma}(x_{\tau}),
\]
where $i^{*}_{\tau,\sigma}$ is the induced map.  Since
$K$ has a unique representative for all $W$-cells in $Y$ it follows that
in fact $i^{*}_{\tau,\sigma}:R(T)^{W_{\tau}} \to R(T)^{W_{\sigma}}$
is the restriction map. Similarly,
\[
C_{W}^{n}(Y;\lZ_{W})=\bigoplus_{\sigma\in I_{n}(K)}
R(G)\otimes \Z[W]^{W_{\sigma}}.
\]
With this description it is clear that
\[
\psi_{n}:=\bigoplus_{\sigma\in I_{n}(K)}\psi_{W_{\sigma}}:
C_{W}^{n}(Y;\lR_{T})\to C_{W}^{n}(Y;\lZ_{W})
\]
is an isomorphism of $R(G)$--modules and the commutativity of
(\ref{commute iso}) shows that this defines an isomorphism of
cochain complexes over $R(G)$. Finally the proof ends by noting
that the cochain complexes $C_{W}^{n}(Y;\lZ_{W})$ and
$C^{n}(Y;R(G))$ are isomorphic and in particular there
is an isomorphism of $R(G)$--modules
$H^{*}_{W}(Y;\lZ_{W})\cong H^{*}(Y;R(G))$.
\qed

\medskip

The hypotheses in the previous theorem can be relaxed if we work
with rational coefficients. More precisely, consider the
coefficient system $\lR_{T}\otimes \Q$ defined by
\[
\lR_{T}\otimes \Q(W/W_{i})=(R(T)\otimes \Q)^{W_{i}}\cong
R(T)^{W_{i}}\otimes \Q.
\]
For this coefficient system we have the following.

\begin{theorem}\label{Bredon cohomology with rational}
Suppose that $G$ is a compact connected Lie group. Let $Y$ be a
$W$--CW complex of finite type. Then there is an isomorphism of
$R(G)\otimes \Q$--modules $H^{*}_{W}(Y;\lR_{T}\otimes \Q)\cong
H^{*}(Y;\Q)\otimes R(G)$.
\end{theorem}
\Proof
Consider the $W$--module $R(T)\otimes \Q$. By definition
\[
\lR_{T}\otimes \Q(W/W_{i})=(R(T)\otimes \Q)^{W_{i}}\cong
R(T)^{W_{i}}\otimes \Q.
\]
For such coefficient systems,
as pointed out in \cite[I. 9]{Bredon}, there is an isomorphism of
cochain complexes
\[
C_{W}^{*}(Y;\lR_{T}\otimes \Q)\cong
\Hom_{\Z[W]}(C_{*}(Y);R(T)\otimes \Q)
\cong \Hom_{\Z[W]}(C_{*}(Y);R(T))\otimes \Q.
\]
It follows that $H_{W}^{*}(Y;\lR_{T}\otimes \Q)\cong
H^{*}(\Hom_{\Z[W]}(C_{*}(Y),R(T))\otimes \Q)$. Let
$D^{*}=\Hom(C_{*}(Y);R(T))$. This cochain complex has a linear
action of $W$ defined by $(w\cdot f)(x)=wf(w^{-1}x)$. Under this
action we have an isomorphism of cochain complexes
\[
\Hom_{\Z[W]}(C_{*}(Y),R(T))\otimes \Q= (D^{*})^{W}\otimes \Q
\cong (D^{*}\otimes \Q)^{W}.
\]
Consider $H^{*}(W;D^{*}\otimes \Q)$; as usual, there are two
spectral sequences computing this group cohomology with coefficients
in a cochain complex. On the one hand, we have
\[
E_{2}^{p,q}=H^{p}(W;H^{q}(D^{*}\otimes \Q))
\Longrightarrow H^{p+q}(W;D^{*}\otimes \Q).
\]
Since we are working with rational coefficients it follows that
\begin{equation*}
E_{2}^{p,q}=H^{p}(W;H^{q}(D^{*}\otimes \Q))\cong \left\{
\begin{array}{ccc}
H^{q}(D^{*}\otimes \Q)^{W}
&\text{ if } &p=0,\\
0& \text{ if } &p>0 .
\end{array}%
\right.
\end{equation*}
On the other hand, we have a spectral sequence
\[
E_{1}^{p,q}=H^{q}(W;D^{p}\otimes \Q)
\Longrightarrow H^{p+q}(W;D^{*}\otimes \Q).
\]
with the differential $d_{1}$ induced by the differential
of the cochain complex $D^{*}$. In this case
\begin{equation*}
H^{q}(W;D^{p}\otimes \Q)\cong \left\{
\begin{array}{ccc}
(D^{p}\otimes \Q)^{W}
&\text{ if } &q=0,\\
0& \text{ if } &q>0.
\end{array}%
\right.
\end{equation*}
Thus the $E_{2}$-term of this spectral sequence is given by
\begin{equation*}
E_{2}^{p,q}= \left\{
\begin{array}{ccc}
H^{p}((D^{*}\otimes \Q)^{W})
&\text{ if } &q=0,\\
0& \text{ if } &q>0 .
\end{array}%
\right.
\end{equation*}
Both of these spectral sequences collapse on the $E_{2}$-term without
extension problems and both converge to $H^{*}(W;D^{*}\otimes \Q)$.
It follows that there is an isomorphism of $R(G)\otimes \Q$ modules
\begin{equation}\label{isom 0}
H_{W}^{*}(Y;\lR_{T}\otimes \Q)=H^{*}((D^{*}\otimes \Q)^{W})
\cong H^{*}(D^{*}\otimes \Q)^{W}=(H^{*}(Y;R(T)\otimes \Q))^{W}.
\end{equation}
The universal coefficient theorem shows that
$H^{*}(Y;R(T)\otimes \Q)\cong H^{*}(Y;\Q)\otimes R(T)$,
with $W$-acting diagonally and thus
\begin{equation}\label{isom 1}
H^{*}(Y;R(T)\otimes \Q)^{W}\cong (H^{*}(Y;\Q)\otimes R(T))^{W}.
\end{equation}
By assumption $Y$ has finite type and thus $H^{*}(Y;\Q)$ is a
finite dimensional $W$--representation over $\Q$.
Using Theorem \ref{W-invariants} below, we obtain
an isomorphism of $R(G)\otimes \Q$--modules
\begin{equation}\label{isom 2}
 (H^{*}(Y;\Q)\otimes R(T))^{W}\cong H^{*}(Y;\Q)\otimes R(G).
\end{equation}
The theorem follows from (\ref{isom 0}), (\ref{isom 1}) and 
(\ref{isom 2}).
\qed

\begin{theorem}\label{W-invariants}
Let $G$ be a compact connected Lie group with torsion--free
fundamental group. Fix $T\subset G$ a maximal torus and let $W$ be
the corresponding Weyl group. If $A$ is any finite dimensional
$W$--representation over $\Q$ then $(A\otimes R(T))^{W}$ is
isomorphic to $A\otimes R(G)$ as modules over $R(G)\otimes \Q$.
\end{theorem}
\Proof Let's consider first the particular case where $A$ is the
regular representation; that is, $A=\Q[W]$.  To prove the theorem in
this case let's consider $R(T)_{0}$ to be the $R(G)$--module $R(T)$
endowed with the \textit{trivial} $W$-action. Define a homomorphism
of $R(G)\otimes \Q$--modules
\begin{align*}
\varphi:\Q[W]\otimes R(T)_{0}&\to \Q[W]\otimes R(T)\\
v\otimes m&\mapsto v\otimes vm,
\end{align*}
for $v\in W$ and $m\in R(T)_{0}$. It is easy to
see that $\varphi$ is a bijection. Moreover, for
any $w,v\in W$ and $m\in R(T)_{0}$ we have
\begin{align*}
\varphi(w\cdot(v\otimes m))&=\varphi((w\cdot v)\otimes m)=
wv\otimes (wv\cdot m)\\
&=w\cdot(v\otimes (v\cdot m))
=w\cdot \varphi(v\otimes m ).
\end{align*}
Therefore $\varphi$ is an isomorphism of $W$--modules over
$R(G)\otimes\Q$. In particular, there is an isomorphism of
$R(G)\otimes\Q$--modules
$$\Q\otimes R(T)_{0} \cong (\Q[W]\otimes
R(T)_{0})^{W}\cong (\Q[W]\otimes R(T))^{W}.$$ Since $R(T)_{0}$ is a
free module over $R(G)$ of rank $|W|$ by \cite[Theorem 1]{Pittie},
we conclude that there is an isomorphism of
$R(G)\otimes\Q$--modules, $(\Q[W]\otimes R(T))^{W}\cong \Q[W]\otimes
R(G)$. This proves the theorem for $A=\Q[W]$. Suppose now that $A$
is any finite dimensional $W$--module over $\Q$. Since $W$ is a
finite group then $A$ is a projective module over $\Q[W]$ and thus
we can find some $\Q[W]$--module $B$ and some integer $m\ge 1$ such
that $A\oplus B\cong (\Q[W])^{m}$. Note that $(A\otimes
R(T))^{W}\oplus (B\otimes R(T))^{W} =((\Q[W])^{m}\otimes R(T))^{W}$
and since the theorem is true for $\Q[W]$, we conclude that
$(A\otimes R(T))^{W}\oplus (B\otimes R(T))^{W}$ is a free module
over $R(G)\otimes \Q$. This means that $(A\otimes R(T))^{W}$ is a
projective module over $R(G)\otimes \Q$. On the other hand, since
$G$ is a compact connected Lie group with $\pi_{1}(G)$ torsion--free
we have
\[
R(G)=\Z[x_{1},\dots, x_{k}]\otimes \Z[y_{1}^{\pm 1},\dots, y_{l}^{\pm 1}]
\]
with $k+l=\mathrm{rank}(G)$ (see for example \cite[Section
1]{Pittie}). We conclude that $R(G)\otimes \Q$ is a tensor product
of a polynomial ring and a Laurent polynomial ring over $\Q$. In
\cite{Swan} Swan extended the Quillen-Suslin theorem to these kind
of rings. Thus any projective  $R(G)\otimes \Q$--module is free. In
particular $(A\otimes R(T))^{W}$ is a free $R(G)\otimes \Q$--module.

\medskip

As a final step we show that $\mathrm{rank}_{R(G)\otimes
\Q}(A\otimes R(T))^{W}=\dim_{\Q}A$. Assuming this is true, then both
$(A\otimes R(T))^{W}$ and $A\otimes R(G)$ are free modules over
$R(G)\otimes \Q$ of
 the same rank and hence isomorphic, proving the theorem.

\medskip

Let us denote by $k$ the fraction field of $R(G)$ (which is is also
the fraction field of $R(G)\otimes \Q$). Because $R(G)$ is an
invariant subring, $K:=R(T)\otimes_{R(G)} k$ is the fraction field
of $R(T)$. Therefore it suffices to prove that $\dim_{k}((A\otimes
k)\otimes_{k} K)^{W}=\dim_{k}(A\otimes k)=\dim_{\Q}A$. To see this
note that  $A\otimes k$ is a finite dimensional $W$--representation
over $k$ and $K$ is a finite Galois extension of $k$ with Galois
group $W$. The result follows applying the lemma given below. \qed

\begin{lemma}
Suppose that $k$ is a field of characteristic $0$. Let $K$ be a
finite Galois extension of $k$ with Galois group $G$.  Assume that
$A$ is a finite dimensional $G$--representation over $k$. If $G$
acts diagonally on $A\otimes_{k} K$ then $\dim_{k}(A\otimes_{k}
K)^{G}=\dim_{k}A$.
\end{lemma}
\Proof Since $k$ has characteristic $0$ then we can find a normal
basis of $K/k$; that is, we can find some $x\in K$ such that the
collection $\{gx\}_{g\in G}$ forms a basis of $K$ as a $k$-vector
space. Using this basis we conclude that the  field $K$, seen as a
$G$--representation over $k$, is such that $K\cong k[G]$. Denote by
$A_{0}$ the vector space $A$ endowed with the \textit{trivial}
$G$-action. Then as in the previous theorem there is an isomorphism
of $G$--representations $A\otimes k[G]\cong A_{0}\otimes k[G]$. In
particular $\dim_{k}(A\otimes K)^{G}=\dim_{k}(A\otimes k[G])^{G}
=\dim_{k}(A_{0}\otimes k[G])^{G}=\dim_{k}A_{0}=\dim_{k}A$.
\qed

\section{Equivariant $K$-theory}\label{equivariant K-theory}

Given a compact $G$--space $X$ the complex equivariant $K$-theory of
$X$, denoted by $K^{0}_{G}(X)$, is defined to be the Grothendieck
group associated to the semi-ring of isomorphism classes of
$G$-equivariant complex vector bundles over $X$. As usual, if $X$ is
a based $G$--space with $G$ acting trivially on the base point
$x_{0}$, the reduced equivariant $K$-theory of $X$ is defined to be
$\tilde{K}^{0}_{G}(X)= \ker(i^{*}:K_{G}^{0}(X)\to
K^{0}_{G}(\{x_{0}\})=R(G))$, where $i:\{x_{0}\}\to X$ is the
inclusion map. When $q\ge 0$ we define
$K_{G}^{q}(X)=\tilde{K}_{G}(\S^{q}\wedge X_{+})$. Equivariant
$K$-theory is representable; that is, for every $q$ there exists a
$G$--space $\bK_{G}^{q}$ such that for every compact $G$--space $X$
there is a natural isomorphism
$\Phi_{X}:K_{G}^{q}(X)\stackrel{\cong}{\rightarrow}[X,\bK_{G}^{q}]_{G}$.
Here $[X,\bK_{G}^{q}]_{G}$ denotes the set of $G$-equivariant
homotopy classes of $G$-equivariant maps $f:X\to\bK_{G}^{q}$. When
$X$ is not a compact space, $K_{G}^{q}(X)$ is defined to be
$[X,\bK_{G}^{q}]_{G}$. For a $G$--space $X$ and $q\ge 0$ there is a
natural isomorphism $K_{G}^{q}(X)\cong K^{q+2}_{G}(X)$ and we denote
$K_{G}^{*}(X)=K^{0}_{G}(X)\oplus K_{G}^{1}(X)$. The projection to a
point $X\to \{x_0\}$ induces a homomorphism of rings
$R(G)=K_{G}^{*}(\{x_0\})\to K_{G}^{*}(X)$ and thus $K_{G}^{*}(X)$
naturally has the structure of a $\Z/2$-graded $R(G)$-algebra, where
$R(G)$ is the representation ring of $G$. Recall that when $G$ is a
compact Lie group $R(G)$ is a domain and also a Noetherian ring (see
for example \cite[Corollary 3.3]{SegalRep}).

\begin{definition}
Suppose that $R$ is a domain. If $M$ is a finitely generated
$R$--module then the rank of $M$, denoted by ${\rm{rank}}_{R}(M)$,
is defined to be largest integer $n$ for which there is an injective
homomorphism of $R$--modules $i:R^{n}\to M$.
\end{definition}

It is easy to see that if $K$ is the field of quotients of $R$ and
$M$ is a finitely generated $R$--module, then the rank of $M$ is the
dimension of $M \otimes K$ as a vector space over $K$. In
particular, if $f:M\to N$ is a homomorphism of finitely generated
$R$--modules such that $f$ induces an isomorphism after passing to
the modules of fractions then
${\rm{rank}}_{R}(M)={\rm{rank}}_{R}(N)$.

\begin{lemma}\label{rankasmodule}
Let $G$ be a compact, connected Lie group with $\pi_{1}(G)$
torsion--free and let $T\subset G$ be a maximal torus. For any
compact $G$--space $X$
${\rm{rank}}_{R(G)}K_{G}^{*}(X)={\rm{rank}}_{\Z}K^{*}(X^{T})$.
\end{lemma}
\Proof Assume that $G$ is a compact Lie group with $\pi_{1}(G)$
torsion--free. The compactness of $X$ implies that $K_{G}^{*}(X)$ is
a finitely generated $R(G)$--module. Also Hodgkin's spectral
sequence gives an isomorphism of $R(T)$--modules (see for example
\cite[Lemma 2.5]{Brylinski}, \cite{Hodgkin})
\begin{equation}\label{passage to T}
K_{T}^{*}(X)\cong K_{G}^{*}(X)\otimes_{R(G)}R(T).
\end{equation}
In particular, since $R(T)$ is a free module over $R(G)$ by
\cite[Theorem 1]{Pittie}, it follows that
\begin{equation}\label{equ1}
{\rm{rank}}_{R(G)}K_{G}^{*}(X)={\rm{rank}}_{R(T)}K_{T}^{*}(X).
\end{equation}
On the other hand, as an application of the localization theorem
\cite[Theorem 4.1]{Segal}, if $i:X^{T}\to X$ is the inclusion map, then
$i^{*}:K_{T}^{*}(X)\to K_{T}^{*}(X^{T})$ is a homomorphism of
$R(T)$-algebras that induces an isomorphism after passing to the modules
of fractions.  Using the previous comment we see
\begin{align}\label{equ2}
{\rm{rank}}_{R(T)}K_{T}^{*}(X)&={\rm{rank}}_{R(T)}K_{T}^{*}(X^{T})\\
&={\rm{rank}}_{R(T)}R(T)\otimes K^{*}(X^{T})
={\rm{rank}}_{\Z}K^{*}(X^{T}).
\end{align}
The lemma follows from equations (\ref{equ1}) and (\ref{equ2}).
\qed

\medskip

In general, if $G$ is a compact Lie group and $X$ is a $G$--CW complex
then associated to the skeletal filtration
\[
X^{0}\subset X^{1}\subset \cdots \subset X^{n}\subset \cdots \subset X
\]
there is a multiplicative spectral sequence \cite{Segal} with
\[
E_{2}^{p,q}=H_{G}^{p}(X;\K^{q}_{G})\Longrightarrow K^{p+q}_{G}(X)
\]
where $\K_{G}^{q}$ denotes the
coefficient system defined by $G/H\mapsto K^{q}_{G}(G/H)$.

\begin{theorem}\label{maintheorem1}
Let $G$ be a compact connected Lie group with $\pi_{1}(G)$
torsion--free and $T\subset G$ a maximal torus. Suppose that $X$ is
a compact $G$--CW complex with connected maximal rank isotropy
subgroups. Assume that there is a CW-subcomplex $K$ of $X^{T}$ such
that for every element $x\in X^{T}$ there is a unique $w\in W$ such
that $wx\in K$ and the family $\{W_{\sigma}~~|~~ \sigma \text{ is a
cell in } K\}$ is contained in a family $\W$ of reflection subgroups
of $W$ satisfying the coset intersection property. If in addition
$H^{*}(X^{T};\Z)$ is torsion--free, then $K_{G}^{*}(X)$ is a free
module over $R(G)$ of rank equal to $\sum_{i\ge
0}{\rm{rank}}_{\Z}H^{i}(X^{T};\Z)$.
\end{theorem}

\Proof Fix $T\subset G$ a maximal torus and let $W$ be the
associated Weyl group. After replacing $X$ with an equivalent
$G$--CW complex, we may assume that the cells of $X$ are of the form
$G/H\times \D^{n}$, where $H$ is a connected subgroup with $T\subset
H$. Consider the spectral sequence associated to the skeletal
filtration of $X$
\[
E_{2}^{p,q}=H_{G}^{p}(X;\K_{G}^{q})\Longrightarrow K_{G}^{p+q}(X)
\]
Under the given hypotheses we have the following:

\noindent{\bf{Claim:}} The $E_{2}$-term of this spectral sequence is
given by
\begin{equation*}
E^{p,q}_{2}\cong \left\{
\begin{array}{ccc}
H^{p}(X^{T};R(G))
&\text{ if } &q \text{ is even},\\
0& \text{ if } &q \text{ is odd}.
\end{array}
\right.
\end{equation*}

Suppose for a moment that the claim is true. The universal
coefficient theorem provides an isomorphism $H^{*}(X^{T};R(G))\cong
H^{*}(X^{T};\Z)\otimes R(G)$. In particular, $H^{*}(X^{T};R(G))$ is
a free module over $R(G)$ of rank $n:=\sum_{i\ge
0}{\rm{rank}}_{\Z}H^{i}(X^{T};\Z)$ and the same is true for
$E^{*,*}_{2}$. The spectral sequence converges to $K_{G}^{*}(X)$ and
by Lemma \ref{rankasmodule},
${\rm{rank}}_{R(G)}K_{G}^{*}(X)={\rm{rank}}_{\Z}K^{*}(X^{T})$. Using
the Chern character we see that ${\rm{rank}}_{\Z}K^{*}(X^{T})
=\sum_{i\ge 0}{\rm{rank}}_{\Z}H^{i}(X^{T};\Z)=n$. Thus
$K_{G}^{*}(X)$ and $E^{*,*}_{2}$ have the same rank as
$R(G)$--modules and $E^{*,*}_{2}$ is a free $R(G)$--module. This
shows that all the differentials $\{d_{r}\}_{r\ge 2}$  must be
trivial and the spectral sequence collapses on the $E_{2}$-term.
Since $E_{2}^{*,*}$ is a free module over $R(G)$ and the spectral
sequence is a sequence of $R(G)$--modules, there are no extension
problems and we conclude that $K_{G}^{*}(X)$ is a free module over
$R(G)$ of rank $n$.

\medskip

We now prove the claim. Suppose first that $q$ is odd; in this case
$\K_{G}^{q}(G/H)=K_{H}^{q}(*)=0$. Therefore when $q$ is odd
$\K_{G}^{q}$ is the trivial coefficient system and in particular
$H_{G}^{p}(X,\K_{G}^{q})=0$. Suppose now that $q$ is even. In this
case $\K_{G}^{q}(G/H)=K^{q}_{H}(*)=R(H)$. This shows that when $q$
is even $H_{G}^{p}(X;\K_{G}^{q})=H_{G}^{p}(X;\lR)$, where $\lR$ is
the coefficient system defined by $\lR(G/H)=R(H)$. Next we show that
there is a natural isomorphism
\begin{equation}\label{step1}
H_{G}^{*}(X;\lR)\cong H^{*}_{W}(X^{T};\lR_{T}).
\end{equation}
To see this note that every cell in $X$ is of the form $G/H\times
\D^{n}$ with $T\subset H$. For such subgroup we have a natural ring
isomorphism $\varphi_{H}:R(H)\to R(T)^{WH}$. These isomorphisms can
be assembled to obtain an isomorphism of cochain complexes
\[
C_{G}^{n}(X;\lR)\cong\bigoplus_{\sigma\in S_{n}(X)} R(G_{\sigma})
\stackrel{\oplus \varphi_{G_{\sigma}}}{\longrightarrow}
\bigoplus_{\sigma\in S_{n}(X)} R(T)^{WG_{\sigma}}\cong
C^{n}_{W}(X^{T};\lR_{T}),
\]
where $S_{n}(X)$ is a chosen set of representatives of all
$n$-dimensional $G$-cells of $X$. In particular, there is an isomorphism
of $R(G)$--modules $H_{G}^{*}(X;\lR)\cong H^{*}_{W}(X^{T};\lR_{T})$.
Finally, note that $X^{T}$ is a compact $W$--CW complex whose $n$-cells
are of the form $W/WH\times \D^{n}$, for some connected maximal rank
isotropy subgroups $H\subset G$. In particular, $WH$ is a reflection
subgroup of $W$ and by  Theorem \ref{Bredon cohomology}
there is an isomorphism of $R(G)$--modules
\begin{equation}\label{step2}
H_{W}^{*}(X^{T};\lR_{T})\cong H^{*}(X^{T};R(G))
\end{equation}
The claim follows from (\ref{step1}) and (\ref{step2}). This
proves the theorem.
\qed

\medskip

With rational coefficients we have a similar spectral sequence
\[
E_{2}^{p,q}=H_{G}^{p}(X;\K^{q}_{G}\otimes \Q)\Longrightarrow
K_{G}^{p+q}(X)\otimes \Q.
\]
In this case the situation simplifies and we have the following theorem.

\begin{theorem}\label{maintheorem2}
Let $G$ be a compact connected Lie group with $\pi_{1}(G)$
torsion--free and $T\subset G$ a maximal torus. Suppose that $X$ is
a compact $G$--CW complex with connected maximal rank isotropy
subgroups. Then $K^{*}_{G}(X)\otimes \Q$ is a free module over
$R(G)\otimes \Q$ of rank equal to $\sum_{i\ge
0}{\rm{rank}}_{\Q}H^{i}(X^{T};\Q)$.
\end{theorem}
\Proof As in the previous theorem we may assume that $X$ is a
$G$--CW complex whose cells are of the form $G/H\times \D^{n}$,
where $H$ is a connected subgroup with $T\subset H$, for a fixed
maximal torus $T\subset G$. Consider the spectral sequence
\[
E_{2}^{p,q}=H_{G}^{p}(X;\K_{G}^{q}\otimes \Q)
\Longrightarrow K_{G}^{p+q}(X)\otimes \Q.
\]
As before we have $E_{2}^{p,q}=0$ when $q$ is odd and when $q$ is
even there is an isomorphism
$E^{p,q}_{2}=H_{G}^{p}(X^{T};\K_{G}^{q}\otimes \Q)\cong
H^{p}_{W}(X^{T};\lR_{T}\otimes \Q)$. Since $X$ is compact then
$X^{T}$ is of finite type. Using Theorem \ref{Bredon cohomology with
rational} we obtain an isomorphism of $R(G)\otimes\Q$--modules,
$H^{*}_{W}(X^{T};\lR_{T}\otimes \Q) \cong H^{*}(X^{T};\Q)\otimes
R(G)$. Thus
\begin{equation*}
E^{p,q}_{2}\cong \left\{
\begin{array}{ccc}
H^{p}(X^{T};\Q)\otimes R(G)
&\text{ if } &q \text{ is even},\\
0& \text{ if } &q \text{ is odd}.
\end{array}
\right.
\end{equation*}
This implies that $E^{*,*}_{2}$ is a free $R(G)\otimes \Q$ module of
rank
\[
n=\sum_{i\ge 0}{\rm{dim}}_{\Q}H^{i}(X^{T};\Q) =\sum_{i\ge
0}{\rm{rank}}_{\Z}H^{i}(X^{T};\Z).
\]
An argument similar to the one provided in Theorem
\ref{maintheorem1} shows that $K_{G}^{*}(X)\otimes \Q$ has rank
equal to $n$ as an $R(G)\otimes\Q$--module, and that the spectral
sequence collapses at the $E_{2}$-term without extension problems.
Therefore $K^{*}_{G}(X)\otimes \Q$ is a free module over
$R(G)\otimes \Q$ of rank $n$. \qed

\section{Applications}\label{applications}

In this section  some applications of Theorems \ref{maintheorem1} and
\ref{maintheorem2} are discussed. We consider examples arising
from representation spheres and the commuting variety in the Lie
algebra of a compact connected Lie group $G$. We also consider
applications arising from inertia spaces, in particular we explore
in detail the case of the space of commuting elements in a Lie
group $G$. Throughout this section $G$ will denote a compact
connected Lie group with torsion--free fundamental group, $T\subset G$
a maximal torus and $W$ the corresponding Weyl group.

\subsection{Conjugation action of $G$ on itself}

Let $G$ act on itself by conjugation. By Example \ref{G on itself}
this action has connected maximal rank isotropy subgroups as we are
assuming that $\pi_{1}(G)$ is torsion--free. In this case $G^{T}=T$
with $W$ acting smoothly on $T$, in particular $T$ has the homotopy
type of a finite $W$--CW complex by \cite[Theorem 1]{Illman}. Note
also that $H^{*}(T;\Z)$ is torsion free and of rank $2^{r}$, where
$r$ is the rank of $G$. Let $\t$ denote the Lie algebra of $T$. Fix
a basis $\Delta$ for the root system $\Phi$ corresponding to the
pair $(G,T)$. Let $\mathfrak{C}(\Delta)$ be the (closed) Weyl
chamber determined by $\Delta$ and $A_{0}$ the unique (closed)
alcove in $\mathfrak{C}(\Delta)$ containing $0\in \t$. Let $K\subset
T$ denote the corresponding alcove in $T$ obtained via the
exponential map. We can give $T$ a $W$--CW complex structure in such
a way that $K$ is CW subcomplex of the underlying CW-complex
structure in $T$. If we further require that $G$ is simply
connected, then as a consequence of \cite[Theorem VII 7.9]{Helgason}
it follows that any element in $T$ has a unique representative in
$K$ under the $W$-action. Also, for every cell $\sigma$ in $K$ we
have $W_{\sigma}=W_{I}$ for some $I\subset \Delta$. Since the family
$\W:=\{W_{I}\}_{I\subset \Delta}$ satisfies the coset intersection
property by Example \ref{example coset}, then it follows that the
hypothesis of Theorem \ref{maintheorem1} are met in this case
yielding the following.

\begin{corollary}\label{Bryl}
Let $G$ be a simply connected compact Lie group act on itself
by conjugation. Then $K_{G}^{*}(G)$ is a free module over
$R(G)$ of rank $2^{r}$, where $r$ denotes the rank of $G$.
\end{corollary}

The previous result was proved in \cite{Brylinski} for the
more general case of compact connected Lie groups with
$\pi_{1}(G)$ torsion--free.

\subsection{Linear Representations}

Let $V$ be a finite dimensional linear representation of $G$
via a homomorphism $\rho:G\to GL(V)$. The action of $G$ on $V$
induces an action of $G$ on $\S^{V}$, the one point compactification
of $V$, with $G$ acting trivially at the point at infinity.
When $V$ is a complex representation the Thom isomorphism
theorem provides an isomorphism of $R(G)$--modules
\begin{equation*}
\tilde{K}^{q}_{G}(\S^{V})\cong \left\{
\begin{array}{ccc}
R(G)&\text{ if } &q \text{ is even},\\
0& \text{ if } &q \text{ is odd}.
\end{array}%
\right.
\end{equation*}

For real representations the coefficient groups
$\tilde{K}_{G}^{*}(\S^{V})$ are not known in general. On the other
hand, we can associate to $V$ another $G$--space by considering the
unit sphere $\S(V)$. To be more precise, endow $V$ with a norm
$\|\cdot\|$ that is invariant under the $G$-action. Such a norm can
be constructed applying the usual trick of averaging any norm with a
suitable Haar measure on $G$. Then $\S(V):=\{x\in V~~|~~ \|x\|=1\}$
is a $G$--space. Note that $\S(1\oplus V)\cong \S^V$. Our main
results imply the following proposition.

\begin{proposition}
Let $V$ denote a representation of $G$ such that the action has
connected, maximal rank isotropy subgroups.
Then $K_G^*(\S(V))\otimes\Q$ and $K_G^*(\S^V)\otimes\Q$ are both free
of rank two as $R(G)\otimes\Q$--modules.
\end{proposition}

Suppose now that $V=\g$, the Lie algebra of $G$ endowed with the
adjoint representation. Given any $X\in \g$ the isotropy subgroup
$G_{X}$ is given by $G_{X}=\{g\in G ~~|~~ Ad_{g}(X)=X\}$. Any
element $X\in \g$ is contained in some Cartan subalgebra $\t\subset
\g$. Let $T$ be the maximal torus in $G$ whose Lie algebra is $\t$
and suppose that $g\in T$. Since $\exp(tX)\in T$ for every $t\in\R$,
then $g$ commutes with $\exp(tX)$ and thus
$\exp(tX)=g\exp(tX)g^{-1}=\exp(Ad_{g}(tX))$ for every $t\in \R$.
This shows that $Ad_{g}(X)=X$ for all $g\in T$. In particular, the
isotropy subgroup $G_{X}$ contains the maximal torus $T$ in $G$.
Also, $G_{X}$ is connected for every $X\in \g$ by \cite[Theorem
3.3.1]{Duistermaat}. This proves that $G$ acts with connected
maximal rank isotropy subgroups on both $\S^{\g}$ and $\S(\g)$. Note
also that $\g^{T}=\t$. Therefore $(\S^{\g})^{T}=\S^{\t}$ and
$\S(\g)^{T}= \S(\t)$ are spheres of dimension $r$ and $r-1$
respectively, where $r$ denotes the rank of $G$. In particular
$H^{*}((\S^{\g})^{T};\Z)$ and $H^{*}((\S(\g))^{T};\Z)$ are both
torsion--free and of rank $2$. Also, these spaces have the structure
of manifolds on which $W$ acts smoothly. By \cite[Theorem 1]{Illman}
it follows that $(\S^{\g})^{T}$ and $\S(\g)^{T}$ have the structure
of a $W$--CW complex. Note that $\t$ can be decomposed into (closed)
Weyl chambers and each (closed) Weyl chamber $\mathfrak{C}(\Delta)$
is determined by a basis $\Delta$ of $\Phi$. For every $x\in
\mathfrak{C}(\Delta)$ the isotropy group $W_{x}$ is a subgroup of
the form $W_{I}$, for some $I\subset \Delta$ and the family
$\W=\{W_{I}\}_{I\subset \Delta}$ satisfies the coset intersection
property by Example \ref{example coset}. Since $W$ acts simply
transitively on the set of all Weyl chambers it follows that we can
obtain $W$--CW complex structures on $(\S^{\g})^{T}$ and
$\S(\g)^{T}$ in such a way that conditions of Theorem
\ref{maintheorem1} are satisfied. As a corollary the following is
obtained.

\begin{corollary}\label{K_{G} lie algebra}
Let $G$ be a compact connected Lie group with $\pi_{1}(G)$ torsion--free.
Suppose that $G$ act on its Lie algebra $\g$ by the adjoint
representation.
If $r$ is the rank of $G$, then
\begin{equation*}
\tilde{K}^{q}_{G}(\S^{\g})\cong \left\{
\begin{array}{ccc}
R(G)&\text{ if } &q\equiv r\ (\text{mod } 2),\\
0& \text{ if } &q+1\equiv r\ (\text{mod } 2).
\end{array}%
\right.
\end{equation*}
Similarly
\begin{equation*}
K^{q}_{G}(\S(\g))\cong \left\{
\begin{array}{cc}
0&\text{ for } q \text{ odd and } r \text{ odd},\\
R(G)&\text{ for } q \text{ odd and }  r \text{ even},\\
R(G)&\text{ for } q \text{ even and }  r \text{ even},\\
R(G)^{2}&\text{ for } q \text{ even and } r \text{ odd}.
\end{array}%
\right.
\end{equation*}
\end{corollary}

\subsection{Arrangements of hyperplanes}

Let $X$ be a $G$--CW complex that satisfies the hypotheses of
Theorem \ref{maintheorem1}. Let $Y\subset X$ be a $G$--CW
subcomplex such that $H^{*}(Y^{T};\Z)$ is torsion--free. Then $Y$
seen as a $G$--space also satisfies the hypotheses of
Theorem \ref{maintheorem1} and in particular $K_{G}^{*}(Y)$ is
a free module over $R(G)$ of rank
$\sum_{j\ge 0}{\rm{rank}}_{\Z}H^{j}(Y^{T};\Z)$. Interesting examples
of $G$--spaces can be obtained in this way. Suppose for example that
$\H=\{H_{a}\}_{a\in \Delta}$ is an arrangement of hyperplanes in $\t$
such that whenever $x\in H_{a}$ for some $a\in \Delta$ and $w\in W$
then $wa\in H_{a'}$ for some $a'\in \Delta$; that is, $\H$ is a
$W$-equivariant arrangement of hyperplanes in $\t$.
Consider the space
\[
C(\H):=\{X\in \t ~~|~~ X\notin \cup_{a\in \Delta} H_{a} \}.
\]
Then $C(\H)$ is a $W$-subspace of $\t$. The connected components of
$C(\H)$ are convex subspaces of $\t$ and thus they are contractible.
In particular, $H^{*}(C(\H);\Z)$ is torsion--free and has rank
$n_{\H}$, where $n_{\H}$ is the number of connected components of
$C(\H)$. Let $X(\H):=\bigcup_{g\in G}Ad_{g}(C(\H))\subset \g$. In
this way we obtain a $G$-subspace of $\g$, with $G$ acting by the
adjoint representation. It is easy to see directly that $X(\H)$ is a
$G$-algebraic variety and in particular it has the homotopy type of
a $G$--CW complex. Note that $X(\H)^{T}=C(\H)$. As a
consequence\footnote{In Theorem \ref{maintheorem1} it is required
that $X^{T}$ is a compact $W$--CW complex. The compactness
hypothesis is used to prove that $X$ is a $G$--CW complex and that
$\sum_{j\ge 0}{\rm{rank}}_{\Z}H^{j}(X^{T};\Z)$ is finite. In this
case both of these conditions are satisfied and thus the theorem can
be applied.} of Theorem \ref{maintheorem1} the following is
obtained.

\begin{corollary}\label{arrangements of hyperplanes}
Let $G$ be a compact connected Lie group with $\pi_{1}(G)$
torsion--free. Let $\H=\{H_{a}\}_{a\in \Delta}$ be a $W$-equivariant
arrangement of hyperplanes in $\t$. Then $K_{G}^{*}(X(\H))$ is a
free $R(G)$--module of rank
$n_{\H}$, where $n_{\H}$ is the number of connected components
of $C(\H)$.
\end{corollary}

\subsection{The commuting variety} Let $G$ be a compact connected
Lie group and $\g$ its Lie algebra. For every integer $n\ge 1$ the
commuting variety in $\g$ is defined to be
\[
C_{n}(\g)=\{(X_{1},\dots,X_{n})\in \g^{n} ~|~ [X_{i},X_{j}]=0
\text{ for all }  1\le i,j \le n \}.
\]
$C_{n}(\g)$ has the structure of an algebraic variety (possibly)
with singularities. The group $G$ acts on $C_{n}(\g)$ via the
diagonal action of the adjoint representation. Let
$(X_{1},\dots,X_{n})\in C_{n}(\g)$ and $g_{i}:=\exp(X_{i})$ for
every $1\le i\le n$. Note that $(g_{1},\dots,g_{n})\in
\Hom(\Z^{n},G)_{\BONE}$. Indeed, the Baker--Campbell--Hausdorff
formula shows that for every $1\le i,j\le n$
\[
g_{i}g_{j}=\exp(X_{i}+X_{j})=\exp(X_{j}+X_{i})=g_{j}g_{i}.
\]

Now the path $\gamma:[0,1]\to \Hom(\Z^{n},G)$ given by $t\mapsto
(\exp(tX_{1}),\dots,\exp(tX_{n}))$ provides a homotopy from the
trivial representation $\BONE$ to $(g_{1},\dots,g_{n})$. In
particular, for every $t\in [0,1]$ the $n$-tuple
$(\exp(tX_{1}),\dots,\exp(tX_{n}))$ is contained in some maximal
torus of $G$ by \cite[Lemma 4.2]{Baird}. Let $\epsilon>0$ be small
enough so that the exponential map is injective on
$B_{\epsilon}(0)$. Choose $0<t<\epsilon$ and let $T$ be a maximal
torus containing $(\exp(tX_{1}),\dots,\exp(tX_{n}))$. Let $\t\subset
\g$ be the Lie algebra of $T$. This is a Cartan subalgebra and since
the exponential map is injective on $B_{\epsilon}(0)$ it follows
that $tX_{i}\in \t$ and thus $X_{i}\in \t$ for all $1\le i\le n$. We
conclude that any $\underline{X}:=(X_{1},\dots,X_{n})\in C_{n}(\g)$
is contained in some Cartan subalgebra in $\g$.\footnote{This
statement is not true for non--compact Lie algebras, for example
this is not true for $\gl_{n}(\C)$.} In particular, the isotropy
subgroup $G_{\underline{X}}$ contains the maximal torus $T$ in $G$
whose Lie algebra is $\t$. Consider $C_{n}(\g)^{+}$, the one point
compactification of $C_{n}(\g)$, with $G$ acting trivially on the
point at infinity. If we assume that $G_{\underline{X}}$ is
connected for every $\underline{X}\in C_{n}(\g)$ then
$C_{n}(\g)^{+}$ is an example of a compact space on which $G$ acts
with connected maximal rank isotropy subgroups. This is the case for
example if $G$ is in the family $\P$.  Indeed, if
$\underline{X}:=(X_{1},\dots,X_{n}) \in C_{n}(\g)$ then
$g\exp(X_{i})g^{-1}=\exp(Ad_{g}(X_{i}))$. Thus $g\in
G_{\underline{X}}$ if and only if $g\in
Z_{G}(\exp(X_{1}),\dots,\exp(X_{n}))$ which is connected and of
maximal rank if $G\in \P$. In this case, given $T\subset G$ a
maximal torus with lie algebra $\t$, then
$(C_{n}(\g)^{+})^{T}=(\t^{n})^{+}=\S^{\t^{n}}$. The Weyl group $W$
acts smoothly on the manifold $\S^{\t^{n}}$ and by \cite[Theorem
1]{Illman} it has the homotopy type of a finite $W$--CW complex. Let
$r$ be the rank of $G$. Then $H^{k}(\S^{\t^{n}};\Q)$ is $0$ unless
$k=0$ or $k=nr$ in which case this is a one dimensional
$W$--representation over $\Q$. This together with Theorem
\ref{maintheorem2} yield the following corollary.

\begin{corollary}
Suppose that $G\in \P$ is a Lie group of rank $r$.
Then there is an isomorphism of $R(G)\otimes \Q$--modules
\begin{equation*}
\tilde{K}^{q}_{G}(C_{n}(\g)^{+})\otimes \Q\cong \left\{
\begin{array}{ccc}
R(G)\otimes \Q&\text{ if } &q\equiv rn\ (\text{mod } 2),\\
0& \text{ if } &q+1\equiv rn\ (\text{mod } 2).
\end{array}
\right.
\end{equation*}
\end{corollary}

\subsection{Inertia spaces}\label{Inertia spaces}
Let $X$ be a compact $G$--CW complex. As we have seen in Section
\ref{section max rank}, if $G$ acts on $X$ with connected maximal
rank isotropy subgroups such that $\pi_1(G_x)$ is torsion--free for
all $x\in X$, then the action of $G$ on the inertia space $\Lambda
X$ also has connected maximal rank isotropy subgroups. Applying
Theorem \ref{maintheorem2} yields the following

\begin{theorem}
Let $X$ denote a compact $G$--CW complex with connected maximal rank
isotropy subgroups all of which have torsion--free fundamental
group. Then $K_G^*(\Lambda X)\otimes \Q$ (as an ungraded module) is
a free $R(G)\otimes \Q$--module of rank equal to $2^{r}\cdot
\left(\sum_{i\ge 0}{\rm{dim}}_{\Q}H^{i}(X^T;\Q)\right)$, where $r$
is the rank of $G$.
\end{theorem}

Suppose now that a compact $G$--CW complex $X$ is such that all of
its isotropy subgroups are of maximal rank and belong to the family
$\P$. Then by Proposition \ref{action on inertia} it follows that
$G$ acts on $\Lambda^{n}(X)$ with connected maximal rank isotropy
subgroups for every $n\ge 0$. Note that if $T\subset G$ is a maximal
torus, then $(\Lambda^{n}(X))^{T}=X^{T}\times T^{n}$. These remarks
together with Theorem \ref{maintheorem2} imply the following result.

\begin{theorem}\label{case of inertia spaces}
Let $X$ denote a compact $G$--CW complex such that all of its
isotropy subgroups lie in $\P$ and are of maximal rank. Then
$K_G^*(\Lambda^n(X))\otimes\Q$ is a free $R(G)\otimes\Q$--module of
rank equal to $2^{nr}\cdot \left(\sum_{i\ge 0}
{\rm{dim}}_{\Q}H^{i}(X^T;\Q)\right)$ where $r$ is the rank of $G$.
\end{theorem}

Suppose that $G\in \P$ and consider the particular case $X=\{x_0\}$.
Then $\Lambda^{n}(X)=\Hom(\Z^{n},G)$ with the conjugation
action of $G$. As a corollary of the previous theorem we obtain
the following.

\begin{corollary}\label{K-theory commuting}
Suppose that $G\in \P$ is of rank $r$. Then
$K_{G}^{*}(\Hom(\Z^{n},G))\otimes \Q$ is free of rank $2^{nr}$ as an
$R(G)\otimes\Q$--module.
\end{corollary}

\medskip

\noindent{\bf{Remark:}} The previous result is no longer true if the
assumption that $\pi_{1}(G)$ is torsion--free is removed.
For example, consider the case $n=1$ and suppose that $G=PSU(3)$
acts on itself by conjugation. Then by
\cite[Proposition 7.4]{Brylinski} there is a $G$-equivariant
line bundle $L$ over $G$ such that $L^{\otimes 3}=1$,
$R(G)/\text{Ann}([L]-1)=\Z$ and
\begin{align*}
K^{0}_{G}(G)&=R(G)^{2}\oplus \Z([L]-1)\oplus \Z([L]^{2}-1),\\
K^{1}_{G}(G)&=R(G)^{2}.
\end{align*}
In particular, $K^{*}_{G}(G)$  and $K^{*}_{G}(G)\otimes \Q$
contain torsion as a module over $R(G)$ resp. $R(G)\otimes \Q$.

\medskip

\begin{example}\label{example commuting tuples in SU(2)} In this
example we compute $K_{G}^{*}(\Hom(\Z^{2},G))$ when $G=SU(2)$. Let
$T\cong \S^{1}$ be the maximal torus  consisting of all diagonal
matrices in $SU(2)$. In this case $W=\Z/2$ acting by
complex conjugation on $T$.  The representation ring
$R(T)=\Z[x_{1},x_{2}]/(x_{1}x_{2}=1)$, the action of
$W$ on $R(T)$ permutes $x_{1}$ and $x_{2}$ and
$R(SU(2))=R(T)^{W}=\Z[\sigma]$, where $\sigma=x_{1}+x_{2}$.
To compute $K_{SU(2)}^{*}(\Hom(\Z^{2},SU(2)))$ we use the spectral
sequence
\begin{equation}\label{spectral case SU(2)}
E_{2}^{p,q}=H_{SU(2)}^{p}(\Hom(\Z^{2},SU(2)),\K_{SU(2)}^{q})
\Longrightarrow K_{SU(2)}^{p+q}(\Hom(\Z^{2},SU(2))).
\end{equation}
As in the proof of Theorem \ref{maintheorem1} we note that
$\K_{SU(2)}^{q}$ is the trivial coefficient system when $q$ is odd
and when $q$ is even we have an isomorphism of $R(SU(2))$--modules
$$H_{SU(2)}^{p}(\Hom(\Z^{2},SU(2)),\K_{SU(2)}^{q})\cong
H_{W}^{p}(T^{2},\lR_{T}).$$
By providing an explicit $W$--CW complex
decomposition to $T^{2}$ it can be proved directly that
\begin{equation*}
H_{W}^{n}(T^{2},\lR_{T})\cong
\left\{
\begin{array}{ccc}
R(SU(2))&\text{ if } &n=0,\\
R(SU(2))^{2}&\text{ if } &n=1,\\
M&\text{ if } &n=2,\\
0&\text{ if } &n>2.\\
\end{array}%
\right.
\end{equation*}
Here $M$ is the $R(SU(2))$--module given by $M=R(SU(2))\oplus
R(SU(2))/ \left<(-\sigma,2)\right>$. It follows that the  spectral
sequence (\ref{spectral case SU(2)}) collapses at the $E_{2}$-term.
Also we conclude that $K^{1}_{SU(2)}(\Hom(\Z^{2},SU(2)))\cong
R(SU(2))^{2}$ and there is a short exact sequence of
$R(SU(2))$--modules
\[
0\to E_{\infty}^{2,0}\cong M\to K_{SU(2)}^{0}(\Hom(\Z^{2},G))
\to E_{\infty}^{0,2}\cong R(SU(2))\to 0.
\]
This sequence splits and thus there are isomorphisms of
$R(SU(2))$--modules
\begin{align*}
K_{SU(2)}^{0}(\Hom(\Z^{2},SU(2)))&\cong R(SU(2))\oplus M,\\
K_{SU(2)}^{1}(\Hom(\Z^{2},SU(2)))&\cong R(SU(2))^{2}.
\end{align*}
It is not hard to show that $K_{SU(2)}^{0}(\Hom(\Z^{2},SU(2))$ is
not a free $R(SU(2))$--module but it becomes free if we invert
$2\in R(SU(2))$.
\end{example}

\medskip

Next we study the relationship between $K^{*}_{G}(\Hom(\Z^{n},G))$
and $K^{*}_{G}(G^{n})$. To start note that in \cite{Brylinski} the
structure of $K_{G}^{*}(G)$ as a $\Z/2$--graded $R(G)$-algebra was
computed for compact connected Lie groups $G$ with $\pi_{1}(G)$
torsion--free. In there an explicit isomorphism of $R(G)$-algebras
$\phi:\Omega^{*}_{R(G)}\to K^{*}_G(G)$ was constructed, where
$\Omega^{*}_{R(G)}$ is the algebra of Grothendieck differentials of
the representation $R(G)$ over $\Z$ (see \cite{Brylinski} for the
definition). Consider now that diagonal action of $G$ on the product
$G^{n}$, where $G$ acts by conjugation on each factor. Since
$\pi_{1}(G)$ is assumed to be torsion--free, Hodgkin's spectral
sequence \cite{Hodgkin} yields
\[
E^{*,*}_{2}=\text{Tor}_{R(G)}^{*,*}(K^{*}_{G}(G),K_{G}^{*}(G))
\Longrightarrow K_{G}^{*}(G^{2}).
\]
$K_{G}^{*}(G)$ is free as a module over $R(G)$, therefore this
spectral sequence collapses on the $E_{2}$-term and there is an
isomorphism of $R(G)$-algebras $K_{G}^{*}(G^{2})\cong
K^{*}_{G}(G)\otimes_{R(G)}K^{*}_{G}(G)$. By induction it follows
that for every $n\ge 1$ there is an isomorphism of $R(G)$-algebras
\begin{equation}\label{K-theory product}
K_{G}^{*}(G^{n})\cong \bigotimes_{R(G)}^{n} K^{*}_{G}(G).
\end{equation}

Suppose now that $X\subset G^{n}$ is a compact connected
$G$-subspace such that $X^{T}=T^{n}$ for some maximal torus
$T\subset G$. The space $\Hom(\Z^{n},G)_{\BONE}$ is an example of
such space. Let $i:X\to G^{n}$ denote the inclusion map and consider
the induced map $i^{*}:K_{G}^{*}(G^{n})\to K^{*}_{G}(X)$. By
assumption the inclusion $i$ induces a homeomorphism
$i^{T}:T^{n}=X^{T}\cong (G^{n})^{T}=T^{n}$. In particular,
$(i^{T})^{*}:K_{T}^{*}((G^{n})^{T})\to K^{*}_{T}(X^{T})$ is an
isomorphism and the localization theorem shows that
$i_{T}^{*}:K_{T}^{*}(G^{n})\to K^{*}_{T}(X)$ is an isomorphism after
inverting the elements in $R(T)-\{0\}$. Since $K_{T}^{*}(G^{n})$ is
free as a module over $R(T)$ and $R(T)$ is a domain this implies
that $i_{T}^{*}$ is injective. By (\ref{passage to T}) there is a
commutative diagram
\[
\begin{CD}
K_{T}^{*}(G^{n})@
>{\cong}>>K_{G}^{*}(G^{n})\otimes_{R(G)} R(T)\\
@V{i_{T}^{*}}VV
@VV{i^{*}\otimes 1}V\\
K_{T}^{*}(X)@>{\cong}>>
K_{G}^{*}(X)\otimes_{R(G)} R(T).\\
\end{CD}
\]
We conclude that $i^{*}:K_{G}^{*}(G^{n})\to K_{G}^{*}(X)$ is also
injective. The above is summarized in the following proposition.

\begin{proposition} Suppose that $G$ is a compact connected Lie group
with $\pi_{1}(G)$ torsion--free. Let $X\subset G^{n}$ be a compact
connected $G$-subspace with $X^{T}=T^{n}$. The inclusion map $i:X
\hookrightarrow G^{n}$ induces an injective homomorphism
$i^{*}:K_{G}^{*}(G^{n})\to K^{*}_{G}(X)$.
\end{proposition}

\noindent Note that in particular, $i^{*}:K_{G}^{*}(G^{n})\to
K^{*}_{G}(\Hom(\Z^{n},G))$ is injective.

In the next proposition we study the map $i^{*}$ for the
particular case of $G=SU(2)$.

\begin{proposition}\label{noncommuting}
Let $i:\Hom(\Z^{2},SU(2))\to SU(2)^{2}$ be the inclusion map.
Then
\[
i^{*}:K^{1}_{SU(2)}(SU(2)\times SU(2))\stackrel{\cong}{\rightarrow}
K_{SU(2)}^{1}(\Hom(\Z^{2},SU(2)))
\]
is an isomorphism and there is a short exact sequence of
$R(SU(2))$--modules
\[
0\to K^{0}_{SU(2)}(SU(2)\times SU(2))\stackrel{i^{*}}{\rightarrow}
K_{SU(2)}^{0}(\Hom(\Z^{2},SU(2)))\to R(\Z/2)\to 0.
\]
\end{proposition}
\Proof
We will show that there is an isomorphism of modules over
$R(SU(2))$
\begin{equation}\label{relative K}
K^{q}_{SU(2)}(SU(2)\times SU(2),\Hom(\Z^{2},SU(2)))\cong
\left\{
\begin{array}{ccc}
0&\text{ if } &q \text{ is even},\\
R(\Z/2)&\text{ if } & q \text{ is odd}.
\end{array}%
\right.
\end{equation}
The proposition follows by considering the long exact sequence in
$K_{SU(2)}^{*}$ associated to the pair
$(SU(2)^{2},\Hom(\Z^{2},SU(2)))$. To prove (\ref{relative K}) note
that $\Hom(\Z^{2},SU(2))$ is a closed $SU(2)$-invariant subspace
of the compact space $SU(2)\times SU(2)$, where $SU(2)$ acts by
conjugation on these spaces. By \cite[Proposition 2.9]{Segal}
there is an isomorphism
\[
K^{q}_{SU(2)}(SU(2)\times SU(2),\Hom(\Z^{2},SU(2)))\cong
\tilde{K}_{SU(2)}^{q}(Y^{+}).
\]
Here $Y:=SU(2)\times SU(2)\setminus \Hom(\Z^{2},SU(2))$ is the space
of non-commuting ordered pairs in $SU(2)$. Consider the commutator
map $\partial:SU(2)\times SU(2)\to SU(2)$. This is a
$SU(2)$-equivariant map. As observed in \cite[Proposition 4.7]{AC}
the restriction of $\partial$, $\partial_{|}:Y\to SU(2)\setminus
\{1\}$ is a locally trivial fiber bundle with fiber
$$F:=\partial^{-1}(-1)=\{(x_{1},x_{2})\in SU(2)\times SU(2) ~|~
[x_{1},x_{2}]=-1\}.$$ Moreover, it is easy to see that this is in
fact a locally trivial $SU(2)$-fiber bundle. The Cayley map provides
a $SU(2)$-equivariant homeomorphism $\psi:SU(2)\setminus \{1\}\to
\su_{2}$, where $\su_{2}$ is the Lie algebra of $SU(2)$ endowed with
the adjoint representation. Thus $Y$ is a locally trivial
$SU(2)$-fiber bundle over $\su_{2}$. Since $\su_{2}$ is
$SU(2)$-contractible it follows that there is a proper
$SU(2)$-equivariant homotopy equivalence $Y\simeq \su_{2}\times F$
and thus there is a $SU(2)$-equivariant homotopy equivalence
$$Y^{+}\simeq(\su_{2}\times F)^{+}=\S^{\su_{2}}\wedge F_{+}.$$
Therefore $\tilde{K}^{q}_{SU(2)}(Y^{+})=\tilde{K}^{q}_{SU(2)}
(\S^{\su_{2}}\wedge F_{+})$. Hodgkin's spectral sequence yields

\begin{equation}\label{ss1}
E^{*,*}_{2}=\text{Tor}_{R(SU(2))}^{*,*}
(\tilde{K}^{*}_{SU(2)}(\S^{\su_{2}}),
\tilde{K}_{SU(2)}^{*}(F_{+}))
\Longrightarrow \tilde{K}_{SU(2)}^{*}(\S^{\su_{2}}\wedge F_{+}).
\end{equation}
Using Corollary \ref{K_{G} lie algebra} we see that
$\tilde{K}_{SU(2)}^{q} (\S^{\su_{2}})$ is $0$ when $q$ is even and
$R(SU(2))$ when $q$ is odd. We now compute
$\tilde{K}^{q}_{G}(F_{+})$. Let $(x_{0},y_{0})$ be any pair of
elements in $SU(2)$ such that $[x_{0},y_{0}]=-1$. As pointed out in
\cite[Lemma 6.2]{ACG2} any element in $F$ is conjugated to
$(x_{0},y_{0})$ and $Z_{SU(2)}(x_{0}, y_{0})= Z(SU(2))\cong\Z/2$.
This means that the conjugation of $SU(2)$ on $F$ is transitive and
there is a $SU(2)$-equivariant homeomorphism $SU(2)/Z(SU(2))\cong
F$, with $SU(2)$ acting on $SU(2)/\Z/2$ by left translation.
Therefore
\begin{equation*}
\tilde{K}^{q}_{SU(2)}(F_{+})=K^{q}_{SU(2)}(F)\cong
K^{q}_{SU(2)}(SU(2)/\Z/2)\cong \left\{
\begin{array}{ccc}
R(\Z/2)&\text{ if } & q \text{ is even},\\
0&\text{ if } &q \text{ is odd}.
\end{array}%
\right.
\end{equation*}
The previous computations show that the spectral sequence (\ref{ss1})
collapses on the $E_{2}$-term and trivially there are no extension
problems proving (\ref{relative K}).
\qed

\section{Appendix}\label{appendix}

In this section we review some basic definitions about $G$--CW complexes
and Bredon cohomology. The reader is referred to \cite{Bredon} and
\cite{Luck} for detailed treatments on these topics. We remark that
throughout this article we work in the category of compactly generated
Hausdorff spaces.

\medskip

Let $G$ be a compact Lie group. Recall that a $G$--space $X$ is a
$G$--CW complex if $X$ is the union of sub $G$--spaces $X^{n}$ such
that $X^{0}$ is the disjoint union of orbits $G/H$ and for every
$n\ge 0$ the space $X^{n+1}$ is obtained from $X^{n}$ by attaching
$G$--cells $G/H\times \D^{n+1}$ along attaching $G$-maps $G/H\times
\S^{n}\to X^{n}$. Given a $G$--space $X$ the set of isotropy
subgroups is $\text{Iso}(X) = \{G_{x} ~|~ x\in X\}$. The orbit type
of $X$ is defined to be $\Fi_{G}(X)=\{(H)~|~ H\in \text{Iso}(X)\}$,
where $(H)$ denotes the conjugacy class of $H$ in $G$.  When
$\Fi_{G}(X)$ is a finite set we say that $X$ has finite orbit type.
Given $H\subset G$ consider the subspaces  $X_{H}:=\{x\in X ~~|~~
G_{x}=H\}$, $X^{(H)}:=\{x\in X ~~|~~ (G_{x})\supset (H)\}$ and
$X^{>(H)}:=\{x\in X ~~|~~ (G_{x})\supsetneq (H)\}$.

\medskip

The following technical conditions on a $G$--space $X$ are very convenient
when dealing with the task of proving that $X$ has the homotopy type of a
$G$--CW complex.

\begin{enumerate}[($G$--CW 1)]
\item $X$ has finite orbit type.
\item $X^{>(H)}\to X^{(H)}$ is a $G$-cofibration for
every $H\in \text{Iso}(X)$; that is, the pair $(X^{(H)},X^{>(H)})$
is a $G$-NDR pair.
\item $X_{H}\to X_{H}/(N_{G}H/H)$ is a numerable principal
$N_{G}H/H$-bundle for every  $H\in \text{Iso}(X)$.
\end{enumerate}

The next theorem provides a criterion to determine that a $G$--space
$X$ has the homotopy type of a $G$--CW complex (see
\cite[Corollary 2.8]{Luck}).

\begin{theorem}\label{Corollary Luck}
Let $G$ be a compact Lie group and $Y$ be a $G$--space satisfying
conditions ($G$--CW 1)-($G$--CW 3). Then $Y$ has the $G$-homotopy type
of a $G$--CW complex $X$ with $\text{Iso}(X) = \text{Iso}(Y)$
if and only if $Y^{H}$ has the homotopy type of a CW complex for
any $H\in \text{Iso}(Y).$
\end{theorem}

Let $\O_{G}$ be the orbit category; that is, $\O_{G}$ is a category with
an object $G/H$ for every (closed) subgroup $H\subset G$ and morphisms
the $G$-equivariant maps $\phi:G/H\to G/K$. Given any such morphism
$\phi:G/H\to G/K$ if $\phi(eH)=gK$, then $gHg^{-1}\subset K$. Thus there
is a morphism $\phi:G/H\to G/K$ in $\O_{G}$ if and only if
$gHg^{-1}\subset K$. Let $R$ be a commutative ring. A coefficient system
in the category of $R$--modules is defined to be a contravariant
functor $M:\O_{G}\to \Mod_{R}$, where $\Mod_{R}$ is the category of
$R$--modules. Suppose that $X$ is a $G$--CW complex,
then there is a coefficient system $\underline{C}_{n}(X)$
defined by
\[
\underline{C}_{n}(X)(G/H)=H_{n}((X^{n})^{H},(X^{n-1})^{H}).
\]
The connecting homomorphism associated to the tripe $((X^{n})^{H},
(X^{n-1})^{H},(X^{n-2})^{H})$ gives rise to a map of coefficient systems
\[
d:\underline{C}_{n}(X)\to \underline{C}_{n-1}(X)
\]
and $d^{2}=0$. Let $\Ca_{G}$ be the category of coefficient systems and
natural transformations between them. Given two coefficient systems $M$
and $N$, denote by $\Hom_{\Ca_{G}}(M,N)$ the abelian group of morphisms
in $\Ca_{G}$ from $M$ to $N$.

\medskip

Let $M$ be a coefficient system in the category of $R$--modules and
$X$ a $G$--CW complex, then there is an associated cochain complex
$(C_{G}^{*}(X;M),\delta)$, where
\[
C_{G}^{n}(X;M)=\Hom_{\Ca_{G}}(\underline{C}_{n}(X),M) \text{ and }
\delta=\Hom_{\Ca_{G}}(d,\text{id}).
\]
The Bredon cohomology of $X$, denoted by $H_{G}^{*}(X;M)$, is defined to
be the cohomology of the cochain complex $(C_{G}^{*}(X;M),\delta)$. Note
that $H_{G}^{*}(X;M)$ has naturally the structure of an $R$--module.
The cochain complex $(C_{G}^{*}(X;M),\delta)$ can be explicitly given in
the following way. For $n\ge 0$
\[
C_{G}^{n}(X;M)\cong\bigoplus_{\sigma\in S_{n}(X)} M(G/G_{\sigma}),
\]
where $S_{n}(X)$ is a chosen set of representatives of all
$n$-dimensional $G$-cells of $X$. It is easy to see that this is
independent of the chosen set of representatives of the $G$-cells in $X$.
With this description, if $x\in C_{G}^{n}(X;M)$ then
$\delta(x)\in C_{G}^{n+1}(X;M)$ is given by
\[
\delta(x)_{\sigma}=\sum_{\tau\in S_{n}(X)}
[\tau:\sigma]M(\tau,\sigma)x_{\tau}.
\]
Here $[\tau:\sigma]$ is as usual the degree of the map induced by the
characteristic map associated to $\sigma$. If $[\tau:\sigma]\ne 0$ then
there is a $G$-map $i_{\tau,\sigma}:G/G_{\sigma}\to G/G_{\tau}$ and
$M(\tau,\sigma)$ is defined as
\[
M(\tau,\sigma):=M(i_{\tau,\sigma}):M(G/G_{\tau})\to M(G/G_{\sigma}).
\]

\section*{Acknowledgments}
The first author is grateful to NSERC for its support and the second
author is grateful to the Pacific Institute for the Mathematical
Sciences for hosting him for the period 2008--2011.


\begin{thebibliography}{99}

\bibitem{AC} A. Adem and F. R. Cohen. Commuting elements and spaces of
homomorphisms.  Math. Ann.  338  (2007), 587-626.

\bibitem{ACG2} A. Adem, F. R. Cohen and J.M. G\'omez. Commuting elements
in central products of special unitary groups. Proceedings of the
Edinburgh Mathematical Society (to appear).

\bibitem{Baird} T. Baird. Cohomology of the space of commuting $n$-tuples
in a compact Lie group, Algebraic and Geometric Topology 7 (2007), 737--754.

\bibitem{Bourbaki} N. Bourbaki. Lie groups and Lie algebras. Chapters 7--9,
Elements of Mathematics (Berlin), Springer-Verlag, Berlin, 2005, xii+434.

\bibitem{Bredon} G. E. Bredon. Equivariant cohomology theories.
Lecture Notes in Mathematics, No. 34, Springer-Verlag, Berlin,
(1967), vi+64 pp.

\bibitem{Brylinski} J. L. Brylinski and B. Zhang. Equivariant $K$-theory of
compact connected Lie groups. $K$-Theory. 20 (2000), no. 1, 23--36.

\bibitem{Dieck1} T. tom Dieck. Transformation Groups. de Gruyter Studies
in Mathematics 8, Walter de Gruyter \& Co., Berlin, (1987), x+312.

\bibitem{Dieck} T. tom Dieck. Transformation Groups and Representation Theory
Lecture Notes in Mathematics 766, Springer-Verlag, Berlin,
(1979), viii+309.

\bibitem{Duistermaat} J.J. Duistermaat and J. A. C. Kolk. Lie groups.
Universitext, Springer-Verlag, Berlin, (2000), viii+344.

\bibitem{Hauschild} V. Hauschild. Compact Lie group actions with isotropy
subgroups of maximal rank. Manuscripta Math., 34 (1981), no. 2-3, 355--379.

\bibitem{Helgason} S. Helgason. Differential geometry, Lie groups, and
symmetric spaces. Graduate Studies in Mathematics,
American Mathematical Society, Providence, RI 2001.

\bibitem{Hodgkin} L. Hodgkin and V. Snaith. Topics in K-Theory, Lecture Notes
in Math 496, Springer-Verlag, New York, 1975.

\bibitem{Illman} S. Illman. Smooth equivariant triangulations of $G$-manifolds
for $G$ a finite group, Math. Ann., 233, (1978),no. 3, 199--220.

\bibitem{KS} V. Kac and A. Smilga. Vacuum structure in supersymmetric
Yang-Mills theories with any gauge group, in The many faces of the superworld.
World Sci. Publ. River Edge, NJ, (1999).

\bibitem{Luck} W. L{\"u}ck. Transformation groups and algebraic {$K$}-theory.
Lecture Notes in Mathematics No. 1408, Springer-Verlag, Berlin,
(1989), xii+443.

\bibitem{Pittie} H. Pittie. Homogeneous vector bundles on homogeneous spaces.
Topology 11, (1972), 199--203.

\bibitem{Segal} G. Segal. Equivariant {$K$}-theory. Inst. Hautes \'Etudes Sci.
Publ. Math. 34,(1968), 129--151.

\bibitem{SegalRep} G. Segal. The representation ring of a compact Lie group.
Inst. Hautes \'Etudes Sci. Publ. Math. 34, (1968), 113--128.

\bibitem{Steinberg} R. Steinberg. On a theorem of Pittie. Topology, 14,
(1975),173--177.

\bibitem{Swan}  R. Swan. Projective modules over {L}aurent polynomial rings.
Trans. Amer. Math. Soc., 237, (1978), 111--120.

\bibitem{Uma} V. Uma. Equivariant $K$-theory of compactifications of algebraic
groups, Transform. Groups, 12, (2007), no. 2, 371--406.

\end{thebibliography}
\end{document}